\input amstex

\documentstyle{amsppt}

\nologo

\magnification=\magstep1
\vcorrection{-0.8cm}

\topmatter
\title
Compatible discrete series
\endtitle
\author
Paola Cellini \\
Pierluigi M\"oseneder Frajria\\
Paolo Papi
\endauthor

\address
Paola Cellini \vskip 0.pt Dipartimento di Scienze
   \vskip 0.pt Universit\`a  di Chieti-Pescara \vskip 0.pt Viale Pindaro 42
\vskip 0.pt 65127 Pescara --- ITALY \vskip 0.pt e-mail:{\rm \
cellini\@sci.unich.it }
\endaddress
\bigskip
\bigskip
\address
Pierluigi M\"oseneder Frajria \vskip 0.pt  Politecnico di
Milano\vskip 0.pt Polo regionale di Como\vskip 0.pt Via Valleggio 11
\vskip 0.pt  22100 Como --- ITALY \vskip 0.pt e-mail:{\rm \
frajria\@mate.polimi.it }
\endaddress
\bigskip
\bigskip
\address Paolo Papi
\vskip 0.pt Dipartimento di Matematica Istituto G. Castelnuovo \vskip 0.pt
Universit\`a di Roma "La Sapienza" \vskip 0.pt Piazzale Aldo Moro 5 \vskip 0.pt
00185 Roma --- ITALY \vskip 0.pt e-mail:{\rm \ papi\@mat.uniroma1.it}
\endaddress
\abstract
Several very interesting
results connecting  the theory of abelian
ideals of Borel subalgebras, some ideas of D.
Peterson relating the previous theory to the
combinatorics of affine Weyl groups, and the theory of
discrete series are stated in  a recent
paper (\cite{Ko2}) by B. Kostant.
\par
In this paper we provide proofs for most of Kostant's results
extending them to  $ad$-nilpotent ideals and
develop one direction of Kostant's
investigation, the compatible discrete series.
\endabstract
\endtopmatter

\TagsOnRight

\def\half{{1\over2}}

\def\cent{\text{\rm cent\,}}
\def\endemo{\qed\enddemo}

\def\({\left(}
\def\){\right)}

\def\stab{\text{\rm Stab}}

\def\a{\alpha}
\def\al{\alpha}
\def\b{\beta}

\def\gg{\gamma}
\def\d{\delta}
\def\th{\theta}
\def\l{\lambda}

\def\t{\tau}

\def\o{\omega^\vee}

\def\p{\Phi}

\def\g{\frak g}
\def\h{\frak h}
\def\i{{\frak i}}

\def\n{\frak n}
\def\bb{\frak b}
\def\ka{\frak k}
\def\pp{\frak p}

\def\ss{\frak s}

\def\R{{\Bbb R}}
\def\ganz{\Bbb Z}
\def\nat{\Bbb N}
\def\ganz{\Bbb Z}

\def\C{\Bbb C}
\def\real{\Bbb R}

\def\D{\Delta}
\def\Dp{\Delta^+}
\def\Dm{\Delta^-}
\def\Da{\widehat\Delta}
\def\Dap{\widehat\Delta^+}

\def\Q{ Q^\vee}

\def\Z{\widetilde Z}
\def\I{\Cal I}
\def\II{\I_{ab}}

\def\DR{\Cal R}

\def\({\left(} \def\){\right)}

\def\Comp{\widetilde{Z}_\tau^{\text{cmpt}}}

\def\lie#1{{\frak #1}}
\document
\heading Introduction \endheading

\medskip
This paper arises from the attempt to understand in detail a recent
paper \cite{Ko2} by B. Kostant, which can be
regarded as an extended research announcement of several very interesting
results connecting (at least) three  topics: the theory of abelian
ideals of Borel subalgebras (which originated
from a much earlier paper by Kostant, \cite{Ko1}), some ideas of D.
Peterson relating the previous theory to the
combinatorics of affine Weyl groups, and, finally, the theory of
discrete series.
\par
In this paper we provide proofs for most of Kostant's results and we
develop one direction of Kostant's
investigation, the compatible discrete series. The paper naturally
divides into three parts that we now describe.\par
Let $\g$ be a simple finite dimensional complex Lie algebra, $G$ the
corresponding connected simply connected Lie
group, $\frak b=\h\oplus \frak n$ be a fixed Borel subalgebra (with
$\h$ Cartan and $\frak n$ nilradical), $\D$ the
root system of $(\g,\h)$, $\Dp$ the positive system induced by the
choice of $\bb$, $W,\,\widehat W,\,\widetilde W$
the finite, affine and extended Weyl groups of $\D$ respectively (for
more details on notation see the list at the
end of the introduction). Let $\I$ denote the set of $ad$-nilpotent
ideals of $\bb$, i.e., the ideals of $\bb$
consisting of $ad$-nilpotent elements.
\par
The first section of the paper is
devoted to build up bijections
$$\widetilde Z\leftrightarrow \Cal W\times
P^\vee/\Q\leftrightarrow\Cal I\times Cent(G)$$
where $\widetilde Z$ is the set of points in $P^\vee$ of the simplex
$\{\sigma\in  \h^*_{\Bbb R} \mid (\sigma,\alpha_i)\leq 1 \text{ for
each  }$ $ i\in \{1,\dots ,n\} \text{ and }
(\sigma,\theta)\geq -2\}$ and $\Cal W$ is a suitable subset of
$\widehat W$ in bijection with $\I$.\par
The previous results generalize to the $ad$-nilpotent case the work
of Kostant and Peterson for the
abelian ideals of $\bb$; their results can be easily obtained
replacing $\widetilde Z$ with
$\widetilde Z_{ab}=\left\{\sigma\in P^\vee \mid
(\sigma,\b)\in\{0,-1,1,-2\}\,\forall\,\b\in\Dp\right\}$.

Let us explain in more detail the meaning of the previous bijections.
To any point $z\in\widetilde Z$ we can associate:
\roster
\item an $ad$-nilpotent ideal $\i_z$;
\item an element $w_z\in\widetilde W$.
\endroster
We define a natural   action of $Cent(G)$ on $\widetilde Z$;  $\i_z$ turns 
out to be  constant on the orbits of this action.
It is well known that $\widetilde W$ is a $Cent(G)$-extension of $\widehat W$. 
In fact we have $\widetilde W= \widehat W\rtimes \Omega$, where $\Omega$ is the 
subgroup, isomorphic to $Cent(G)$, of the elements in $\widetilde W$ 
which stabilize fundamental alcove. 
Then $z\mapsto w_z$  maps each $Cent(G)$-orbit of $\widetilde Z$ onto a 
left coset of $\widetilde W/\Omega$. Moreover, if $\widehat w_z$ is the unique 
element in $(w_z\Omega)\cap \widehat W$, then $\i_z\mapsto \widehat w_z$  realizes
a bijection $\I\to\Cal W$. 
\par
Restricting to $z\in \widetilde Z_{ab}$, if we set $v_z^{-1}$ to be the $W$-component of
$w_z\in \widetilde W= W\ltimes T(P^\vee)$, we have that $\i_z$ is the sum of the root spaces $\g_\a$ for all $\a\in
-v_z^{-1}(\D^{-2}_z)\cup v_z^{-1}(\D^{1}_z)$, where by definition
$\D^{i}_z=\{\a\in\Dp\mid (\a,z)=i\}.$

\medskip The
connection of these results with representation theory is the main
theme of section two. Set
$X=\overline C_2\cap P^\vee$
and $dom:\widetilde Z_{ab}\to X$ to be the map defined by
$dom(z)=v_{z}^{-1}(z)$.
Given $\tau\in X$, then
$\Theta_{\tau}=Ad(\exp(\sqrt{-1}\pi\tau)$ is a Cartan involution for
$G$ and the corresponding  Cartan decomposition
$\g=\ka_{\tau}\oplus\lie p_{\tau}$ has the property that $\ka_{\tau}$
is an equal rank symmetric Lie subalgebra of $G$. Given a Cartan
involution $\Theta$ and a $\Theta$-stable Borel subalgebra $\bb'$, we  introduce, after Kostant,  a notion of compatibility of
$\bb'$ with $\Theta$ and we prove  that $\bb'$ is compatible with $\Theta$ if and only if the pair
$(\Theta, \bb')$ is conjugate under $G$ to a pair
$(\Theta_{\tau},\bb)$ with $\tau\in X$.
\par
Now, for $\tau\in X$, set $\widetilde{Z}_\tau=dom^{-1}(\tau)\cap
\widetilde{Z}_{ab}$; we provide proofs for the two important results
of \cite{Ko2}. First, $\widetilde{Z}_\tau$ is in canonical bijection with
$W_\tau\backslash W$, $W_\tau$ being the Weyl group of $(\ka_\tau,\h)$.
Moreover we consider
$$\widetilde{Z}_\tau^{\text{cmpt}}=\{z\in\widetilde{Z}_\tau\mid
v_z^{-1}(\bb) \hbox{ is
compatible with }\Theta_\tau\}$$ and we prove that it has $|Cent(G)|$ elements
by exhibiting a bijection
with a suitable copy of $Cent(G)$ inside $\widetilde W$.
\par
To introduce the connection with discrete series, fix again $\tau\in X$
and set
$G_\tau$ to be the real form of
$G$ corresponding to $\Theta_\tau$. Set  $K$ to be a maximal compact
subgroup of $G_\tau$.
If $\mu$ is a regular integral weight denote by  $\pi_\mu$ 
the corresponding discrete series representation for $G_{\tau}$
in Harish Chandra parametrization.
\par
The previous results imply that, if $\l$ is regular and dominant, then
the map $z\mapsto \pi_{v_{z}^{-1}(\l)}$ defines a bijection between
$\widetilde Z_{\tau}$ and the discrete series having infinitesimal
character $\chi_{\l}$. What is more important, the minimal $K$-type of
$\pi_{v_{z}^{-1}(\l)}$
and the $L^2$-cohomological degree on $G_\tau/T$ in which
$\pi_{v_{z}^{-1}(\l)}$ appears, can be expressed in terms of
combinatorial
data related to $\i_z$.
\par
We can also single out the discrete series
$\pi_{v_{z}^{-1}(\l)}$ that correspond to   $z\in \Comp$. These are
called {\it compatible} discrete series: we note that they are
exactly the small discrete series of Gross and Wallach (\cite{GW}).
We then go on to show that the $K$-spectrum of the compatible discrete
series can be computed fairly explicitly. The relevant results in 
this direction
are given in Theorem~2.8 and the discussion thereafter.
\medskip

Section 3 is devoted to the proofs of the results of \cite{Ko2} that
did not fit in our previous discussion, namely, we develop Kostant's
theory of special and nilradical abelian ideals.
   The main outcome is that the
structure theory for an abelian ideal $\i_z$ can be given a symmetric
space significance in terms of $\tau=dom(z)$
and of the
associate decomposition
$\g=\ka_\tau\oplus\frak p_\tau$. 
More precisely,
if $\i_z=\left(\bigoplus\limits_{\a\in -v_z^{-1}(\D^{-2}_z)}\g_\a\right)\oplus
\left(\bigoplus\limits_{\a\in v_z^{-1}(\D^{1}_z)}\g_\a\right)$, then
the two summands in the r.h.s. of the
previous expression correspond to $\i_z\cap \ka_\tau,\,\i_z\cap\frak
p_\tau$, respectively. Other features of this decomposition are the 
following. On one hand
$\i_z\cap \ka_\tau$ is again an abelian ideal of
$\bb$, and the abelian ideals arising in this way (called {\it
special}) can be abstractly characterized. On the other hand,
$\i_z\cap\frak p_\tau$ is an abelian $\bb\cap\ka_\tau$-submodule of 
$\bb\cap\frak p_\tau$
and the map $z\mapsto \i_z\cap\frak p_\tau$ sets up a bijection between
$\widetilde Z_\tau$ and the set of
such submodules; such a map extends to a bijection between $dom^{-1}(\tau)\cap \widetilde
Z$ and all the 
$\bb\cap\ka_\tau$-submodules of
$\bb\cap\frak p_\tau$. Finally, we give a
criterion to decide whether  a special abelian ideal is the
nilradical of a parabolic subalgebra of $\g$ and we determine
explicitly all such ideals.  All these results are obtained
combining  the results of section 1 with
elementary combinatorics of root systems.
\par\newpage
\heading Notation\endheading
$$\alignat2
&\g\quad &&\text{finite-dimensional complex simple Lie algebra,}\\
&G\quad &&\text{connected  simply connected semisimple Lie}\\
& &&\text{group with Lie algebra $\g$,}\\
&Cent(G)\quad &&\text{center of $G$,}\\
&T\quad &&\text{maximal compact torus of $G$, with Lie algebra $\lie t$,}\\
& &&\text{hence $\h=\lie t\oplus \sqrt{-1}\lie t$ is a Cartan 
subalgebra of $\g$,}\\
&\h_{\Bbb R}=\sqrt{-1}\lie t\\
&\bb=\h\oplus\frak n\quad &&\text{Borel subalgebra of $\g$, with
Cartan component $\h$}\\
& &&\text{and nilradical $\frak n$,}\\
&\D\subset \h_{\Bbb R}^*\quad &&\text{finite (irreducible)
root system of $\g$}\\
& &&\text{with
positive system
$\Dp$,}
\\ &\Pi=\{\a_1,\ldots,\a_n\}\quad &&\text{simple roots of $\Dp$,}\\
&Q=\sum_{i=1}^n\Bbb Z\a_i,\,\Q=\sum_{i=1}^n\Bbb
Z\a_i^\vee&&\text{root and coroot lattices,}\\
&P=\sum_{i=1}^n\Bbb
Z\omega_i,\,P^\vee=\sum_{i=1}^n\Bbb Z\omega_i^\vee&&\text{weight and
coweight lattices,}\\
&\th=\sum_{i=1}^nm_i\a_i\quad &&\text{highest root of $\D$ w.r.t. $\Pi$,}\\
&J=\{i\mid 1\leq i\leq n,\, m_i=1\}\\
&o_i=\omega_i^\vee/m_i,\,1\leq i\leq n&&\text{nonzero vertices of the
fundamental alcove,}\\
&\rho=\omega_1+\ldots+\omega_n&&\text{sum of fundamental weights,}\\
&\rho^\vee=\omega_1^\vee+\ldots+\omega_n^\vee&&\text{sum of
fundamental coweights,}\\
&h\quad &&\text{Coxeter number of $\D$,}\\
&\d\quad &&\text{fundamental imaginary root,}\\
&\Da=\D+\ganz\d\quad &&\text{affine root system associated to $\D$,}\\
&\Dap=(\Dp+\nat \d)\cup
(\Dm+\nat^+\d)\ \
&&\text{positive system for $\Da$,}\\
&W\quad &&\text{finite Weyl group,} \\
&T(L)\subset Aff(\h^*_{\Bbb R})\quad &&\text{group of translations corresponding}\\
& &&\text{to the lattice $L\subset\h^*_{\Bbb R}$,} \\
&\widehat W=W\ltimes T(\Q)\quad &&\text{affine Weyl group,} \\
&\widehat W_i=W\ltimes T(i\Q)\quad &&i\in\nat, i\geq 1,\ (\widehat 
W_1\equiv\widehat W),\\
&\widetilde W=W\ltimes T(P^\vee)\quad &&\text{extended affine Weyl group,} \\
&\I=\{\i\mid \i\text{ ideal of $\bb$, }\i\subset\frak n\}\quad
&&\text{$ad$-nilpotent ideals of
$\bb$,}\\
&\II\quad &&\text{abelian ideals of $\bb$.}
\endalignat$$
\par\noindent
We refer to \cite{CP} for a brief description of $\widehat \D$ and
$\widehat W$, though we do not completely follow the notation we used there.
In particular, in \cite{CP, Section 1} we  described in great detail 
the relationship
between $\widehat W$ viewed as the Weyl group of $\widehat \D$ and its affine
representation on $\h_{\Bbb R}^*$. By virtue of that description, we do not
distinguish here between $\widehat W$ and its affine representation.
We note that $\widehat W_i$ can be viewed as the affine Weyl group of 
${1\over i} \D$.
\smallskip
We will denote by $(\ ,\ )$ at the same time the Killing form on 
$\h_{\Bbb R}$, the
induced form on
$\h_{\Bbb R}^*$ and the natural pairing $\h_{\Bbb R}\times\h_{\Bbb 
R}^*\to \Bbb R$. We set:
$$\alignat2
&H_{\a,k}=\{x\in \h_{\Bbb R}^*\mid (x,\a)=k\},\ \a\in\Dp,\ k\in \ganz\quad
&& \text{affine reflection hyperplane of $\widehat W$,}\\
&C_\infty=\{x\in \h_{\Bbb R}^*\mid (x,\a_i)>0\, \forall\ i=1,\dots ,n\}
\quad &&\text{the fundamental chamber of $W$,}\\
&C_i=\{x\in \h_{\Bbb R}^*\mid (x,\a)>0\,\forall\ \a\in\Pi,\ (x,\theta)<i\}
\quad &&\text{the  fundamental alcove  of $\widehat W_i$.}\\
\endalignat
$$
Recall that  $\D$  can be given a partial order defined by $\a<
\b$ if $\b-\a$ is a sum of positive roots.
We will denote by $V_\a$ the principal dual order ideal generated by
$\a$, namely:
$$V_\a=\{\b\in\Dp\mid\b\geq \a\}.$$
\smallskip
Moreover, for $\i\in\I$ we define $\p_\i=\{\a\in\Dp\mid
\g_\a\subset\i\}$, and for $\Phi\subset \Dp$ we set
$\i(\p)=\bigoplus\limits_{\a\in\p}\g_\a$, so that $\i\mapsto\p_\i,\
\p\mapsto\i(\p)$ are mutually inverse maps between $\I$ and the set of dual order ideals of the poset
$(\Dp,\leq)$. 
We shall  use the notation $\i(\p)$  even when $\p$ is any subset  of $\Dp$ (so $\i(\p)$ is not
necessarily an ideal).\par  We set $\i(\emptyset)$ to be the zero ideal.
\heading\S1 Connections between $ad$-nilpotent ideals, the affine
Weyl group and the center of
$G$\endheading
\bigskip

In this section we will explain how  $ad$-nilpotent ideals can be
naturally   parametri\-zed by a
certain subset $\Cal W$ of $\widehat W$, and in turn, by the points
in $\Q$ of a suitable simplex
$D\subset \h_{\Bbb R}^*$;  on this simplex there is a natural action of
$Cent(G)$, which should be regarded as
embedded into $\widetilde W$: we will describe the relationships of
this action with  the
$ad$-nilpotent ideals. These constructions  specialize nicely to the
case of abelian ideals.
The idea of  relating abelian ideals with elements in  $\widehat W$
appears in \cite{Ko2} (where it is
attributed to D. Peterson); the generalization  to the $ad$-nilpotent
case can be found in
\cite{CP}. The role of the center has been pointed out in \cite{Ko2}
(in the abelian case).
\smallskip

First we explain how to view  $\I$ as a subset of $\widehat W$.

Let
$\i=\bigoplus\limits_{\a\in\p}\g_\a,\,\p\subset\Dp$, be an
$ad$-nilpotent ideal. Set
$$L_\i=\bigcup_{k\geq 1}\left( -\p^k+k\d\right)\subset \Dap,$$
where $\p^1=\p$ and $\p^k=(\p^{k-1}+\p)\cap \D,\,k\geq 2$. Note that,
since $\i$ is nilpotent, $L_\i$ is a finite set. For
$w\in
\widehat W$ set
$$N(w)=\left\{\a\in\Dap\mid w^{-1}(\a)<0 \right\}.$$
Then:

\proclaim{Proposition A} \cite{CP, Theorem 2.6}
There exists a unique  $w_\i\in\widehat W$ such that $L_\i=N(w_\i)$.
\endproclaim

\smallskip
Hence we have an injective map $f:\I\to \widehat W,\,f(\i)=w_\i$. We set $$\Cal W=f(\I), \quad \text{ and } \quad \Cal
W_{ab}=f(\I_{ab}).$$
$\Cal W$ and $\Cal W_{ab}$ are characterized inside $\widehat W$ by 
the following
properties.

\smallskip
\proclaim
{Proposition B} \cite{CP, Theorem 2.9, Proposition 2.12}
\par\noindent{\rm (1)}
We have $$\Cal W_{ab}=\{w\in \widehat W\mid w(C_1)\subset C_2\}.$$
\par\noindent{\rm (2)}
Assume $w\in \widehat W$ and $w=t_\tau v$, with $\tau\in \Q$, $v\in W$.  Then
$w\in\Cal W $  if and only if the following
conditions hold:
\roster
\item"(i)"
$w(C_1)\subset C_{\infty}$;
\item"(ii)"
$(v^{-1}(\tau),\a_i) \leq 1$ for each $i\in\{1,\dots,n\}$ and
$(v^{-1}(\tau), \theta) \geq - 2$.
\endroster
\endproclaim

\smallskip
\proclaim
{Corollary}(D. Peterson) $|\II|=2^n$.
\endproclaim
\medskip

For $\i\in \I$ and $\a\in \Dp$ we have that $\a\in \Phi_\i$
if and only if $-\a+\d\in N(w_\i)$.  As shown in \cite{CP, Section 
1}, we have that
$-\a+\d\in N(w_\i)$ if and only $H_{\a, 1}$ separates $C_1$ and 
$w_\i(C_1)$. We shall need the following result \cite{IM, \S 1.9}.

\smallskip
\proclaim
{Lemma C}
Let  $w\in \widetilde W$, $w=t_\tau v$, $\tau\in P^\vee$, $v\in W$. 
Set $\D_\tau^i=\{\a\in \Dp\mid (\a, \tau)=i\}$. 
Then $H_{\a, 1}$ separates $C_1$ and $w(C_1)$ if and only if 
$\a \in \left(\bigcup\limits_{i>1}\D_\tau^i\right)\cup(\D_\tau^1\setminus 
N(v))$.
\endproclaim

\smallskip
We shall also  need the following classical results, which  can be extracted
from \cite{IM, \S1}, on the embedding of $Cent(G)$ into $W$. Let
$\D(j)$ denote the root subsystem of $\D$ generated by
$\Pi\setminus\{\a_j\}$ and by $w^j_0$ the longest (w.r.t.
$\Pi\setminus\{\a_j\}$) element of the
corresponding parabolic subgroup of $W$. Let $w_{0}$ be the longest
element of $W$ with respect to $\Pi$, and
$$\Omega=\left\{t_{\omega_j^\vee}w^j_0w_0\mid j\in
J\right\}\cup\{1\}\subset\widetilde W.$$

\proclaim{Proposition D}
\par\noindent{\rm (1)}
$\Omega$ is a group. Precisely, $\Omega$ is the subgroup of all elements
$w\in \widetilde W$ such that $w(\overline C_1)=\overline C_1$, hence
it is isomorphic to
$P^\vee/\Q$. Composing this isomorphism with the one induced by the exponential map we also obtain an isomorphism 
$\Omega\to Cent(G)$.
\par\noindent{\rm (2)}
The projection $t_{\omega_j^\vee}w^j_0w_0\mapsto w^j_0w_0$
from $\Omega$ to $W$ is injective, hence  it embeds $Cent(G)$ into $W$.
\endproclaim

\remark{\sl Remark} It is clear that, for any $1\leq j\leq n$,
$N(w_0^jw_0)$ equals  $\Dp\setminus \D(j)^+$. In particular, for 
$j\in J$ we have
$N(w_0^jw_0)=\{\a\in\Dp\mid (\a,\omega_j^\vee)=1\}$.\endremark

\bigskip
The last tool which will be relevant in what follows is  the
encoding of $ad$-nilpotent
ideals by means of lattice points of a certain simplex $D$ in $\h_{\Bbb R}^*$.
The definition of this simplex is motivated by  Proposition B; set
$$D=\left\{\sigma\in  \h_{\Bbb R}^* \mid (\sigma,\alpha_i)\leq 1 
\text{ for each
} i\in \{1,\dots ,n\} \text{ and }
(\sigma,\theta)\geq -2\right\}$$
and define
$$\align
&\widetilde Z=D\cap P^\vee,\\
&Z=D\cap \Q,\\
&\widetilde Z_{ab}=\left\{\sigma\in P^\vee \mid
(\sigma,\b)\in\{0,-1,1,-2\}\,\forall\,\b\in\Dp\right\},\\
&Z_{ab}=\left\{\sigma\in \Q\mid
(\sigma,\b)\in\{0,-1,1,-2\}\,\forall\,\b\in\Dp\right\}.\endalign$$

\medskip
It is easily seen that $D=t_{\rho^\vee}\,w_0(\overline
C_{h+1})=\rho^\vee-\overline C_{h+1}$, $h$ being the Coxeter number of $\D$. Note that
$t_{\rho^\vee}\,w_0\in\widetilde W$; on the other hand in \cite{CP2, Lemma 1}
it is proved that there exists an element $\widetilde w\in\widehat W$
such that $D= \widetilde w\,(\overline C_{h+1})$. Set, for $r\geq 1$
$$\Omega_{r}=\left\{t_{r\omega_j^\vee}w_0^jw_0\mid j\in J\right\}\cup\{1\};$$
then $\Omega_{r}$ is the subgroup of all elements
$w\in \widetilde W$ such that $w(\overline C_{r})=\overline C_{r}$.
It follows that $\widetilde w \,\Omega_{h+1}\, {\widetilde w}^{-1}$
is the subgroup of all elements
$w\in \widetilde W$ such that $w(D)=D$.
\par
We set
$\Sigma=\left\{t_{-\omega_j^\vee}w^j_0w_0\mid j\in J\right\}\cup\{1\}.$
Clearly, $\Sigma$ is a subgroup of $\widetilde W$ isomorphic to $\Omega$.
We have the following result.
\medskip

\proclaim{Proposition 1.1}
\par\noindent{\rm(1)}
$\Sigma=\left\{w\in \widetilde W\mid w(D)=D\right\}$,
so that $\Sigma=\widetilde w\, \Omega_{h+1} \,{\widetilde w}^{-1}$.
In particular $\Sigma$ acts on $\widetilde Z$.
\par\noindent{\rm(2)}
For any $z\in \widetilde Z$ the orbit of $z$ under the action of $\Sigma$ is
a set of representatives of $P^\vee/ Q^\vee$. In particular $\Sigma$
acts freely on  $\widetilde Z$.
\par\noindent{\rm(3)}
The action of $\Sigma$ on $\Z$ preserves $\Z_{ab}$.
\endproclaim

\demo{Proof}
\par\noindent{(1)}  It suffices to prove that any element in $\Sigma$
preserves $D$. Take
$x\in D$ and consider $t_{-\omega_j^\vee}w^j_0w_0(x)$, $j\in J$. For
$i=1,\ldots,n$, we have
$(t_{-\omega_j^\vee}w^j_0w_0(x),\a_i)=(x,w_0w^j_0(\a_i))-(\omega_j^\vee,\a_i)$.
Now if $\a_i\in\D(j)$ we have $(\omega_j^\vee,\a_i)=0$ and that
$w_0w^j_0(\a_i)$ is a simple root (see \cite{IM}). Since $x\in D$, we get
$(x,w_0w^j_0(\a_i))\leq 1$ as desired. If $\a_i\notin\D(j)$, then
$\a_i=\a_j$, so
that $(\omega_j^\vee,\a_i)=1$ and $w_0w^j_0(\a_i)=-\th$; therefore
$(t_{-\omega_j^\vee}w^j_0w_0(x),\a_i)=-(x,\th)-1\leq 2-1=1$. Finally, we have
$$(t_{-\omega_j^\vee}w^j_0w_0(x),\th)=(x,w_0w^j_0(\th))-(\omega_j^\vee,\th)=
(x,w_0(\a_j))-1\geq-1-1=-2.$$ The fact that $\Sigma$ also preserves $\widetilde
Z$ is immediate.
\par\noindent{(2)} Let $z\in \widetilde Z$. For any root $\alpha$, denote
by $s_\alpha$
the reflection associated to $\alpha$; we have
$s_\alpha(z)=z-(z,\alpha^\vee)\alpha=
z-(z,\alpha)\alpha^\vee$ and, since $z\in P^\vee$,
$(z,\alpha)\alpha^\vee\in Q^\vee$.  It
follows that, for any $v\in W$, $z$ and $v(z)$ differ by an element
in $Q^\vee$.
Therefore, for any $j\in J$, $z$ and
$t_{-\omega_j^\vee}w^j_0w_0(z)=-\omega^\vee_j+w^j_0w_0(z)$
are distinct $\mod Q^\vee$, since $\omega^\vee_j\not\in Q^\vee$. The
claim follows directly.
\par\noindent{(3)} If $x\in\widetilde Z_{ab}$, we have $$\multline
(t_{-\omega_i^\vee}w^i_0w_0(x),\b)=(x,w_0w^i_0(\b))-(\omega_i^\vee,\b)\\=
\cases
(x,w_0w^i_0(\b))\in\{0,1,-1,-2\}\quad&\text{if }\b\in\D(i)\\
(x,w_0w^i_0(\b))-1\in
\{0,1,-1,2\}-1=\{0,1,-1,-2\}\quad&\text{otherwise.}\endcases\endmultline$$
Hence
$t_{-\omega_i^\vee}w^i_0w_0(x)\in \Z_{ab}$, as desired. \endemo

\proclaim
{Lemma 1.2} $W\cdot\overline C_{k}=\{x\in \h_{\Bbb R}^*\mid -k\leq (x,\b) \leq
k\ \forall\,\b\in\Dp\}$.
\endproclaim

\demo{Proof}
Consider $x\in W\cdot\overline C_{k}$; then
$x=\sum_{i=1}^n\lambda_i w(o_i)$ for some $w\in W$ and $\lambda_i\geq
0,\,\sum_{i=1}^n\lambda_i\leq k$. Set $\b_i=w(\a_i)$ for $1\leq i\leq n$; then
$\b_1,\ldots,\b_n$ is another basis for $\D$, hence for any $\b\in \Dp$ we have
$\b=\sum_{i=1}^n\mu_i\b_i$ with $\vert \mu_i\vert\leq m_i$ and  $\mu_i\geq 0\
\forall\,i$ or $\mu_i\leq 0\ \forall\,i$. Finally we have
$\vert(x,\b)\vert=\vert \sum\limits_{i=1}^n\lambda_i\frac{\mu_i}{m_i}\vert\leq
\sum_{i=1}^n\lambda_i\leq k$.
\endemo

\proclaim
{Proposition 1.3} For any $z\in\Z$ we have $z+\overline C_1\subset
W\cdot\overline C_h$. Moreover $z\in \Z_{ab}$ if and only if
$z+\overline C_1\subset
W\cdot\overline C_2$. In particular, for any $z\in\Z$ there exists a unique
$v\in W$  such that $v(z+\overline C_1)\subset \overline C_h$ and,
for such a $v$,
$v(z+\overline C_1)\subset \overline C_2$ if and only $z\in \Z_{ab}$.
\endproclaim

\demo{Proof}
Let $z\in \Z$: we may
write $z=\rho^\vee-\sum_{i=1}^n\lambda_io_i,\
\lambda_i\geq0,\,\sum_{i=1}^n\lambda_i\leq h+1$. For any positive root
$\b=\sum_{i=1}^n\mu_i\a_i,\,0\leq \mu_i\leq m_i$ we have
$$(z,\b)=\sum_{i=1}^n\mu_i-\sum_{i=1}^n\lambda_i\frac{\mu_i}{m_i}\leq
\sum_{i=1}^nm_i=h-1.$$
On the other hand, since at least one of the $\mu_i$ is positive,
$(z,\b)\geq 1 -\sum_{i=1}^n\lambda_i\geq 1-(h+1)=-h$.
\par Now  for any $y\in \overline C_1$ and any $\b\in \Dp$ we have
$0\leq (y,\b)\leq
1$, whence $-h\leq (z+y,\b)\leq h$. By  Lemma 1.2 we get $z+\overline
C_1\subset
W\cdot\overline C_h$.
\par
We see directly that  $z\in \Z_{ab}$ if and only if $-2 \leq(z+y,\b)\leq 2$,
or equivalently,  by Lemma 1.2, if and only if $z+\overline
C_1\subset W\cdot\overline C_2$.
\endemo

\medskip\noindent
{\bf Definition.}
Let $z\in \Z$ and $v\in W$ be (the unique element of $W$)
such that $v(z+\overline C_1)\subseteq \overline C_h$. We define the map 
$$\widetilde F: \widetilde Z\to \widetilde W, \quad z\mapsto t_{v(z)} v,
\qquad
\text { and we set } \qquad v_z=v^{-1}. $$
There exists a unique $w\in \widehat  W$ such that $v(z+\overline C_1)
=w(\overline C_1)$, thus we also have a map
$$F:\widetilde Z\to \widehat W,\quad z\mapsto w.$$
(In the Introduction $\widetilde F(z),\,F(z)$ have been called $w_z,\,\widehat w_z$, respectively).\par
By Proposition A we can thus associate to $z$ an ad-nilpotent ideal 
$$\i_{z}=f^{-1}(F(z)).$$

\bigskip
\proclaim{Proposition 1.4}
\par\noindent{(a)} $F(z)=F(z'),\ z,z'\in\Z$  if and only if $z=\psi(z')$ for
$\psi\in\Sigma.$
\newline
\par\noindent{(b)} $F$ is a surjection $\widetilde Z\twoheadrightarrow \Cal W$
and $F_{|Z}:Z\to\Cal W$ is a bijection, with inverse map
$t_{\tau}v\mapsto v^ {-1}(\tau)$. \newline
\par\noindent{(c)} Set $H(z)=(F(z),z\,\,mod\,\Q)$. Then $H$ is
bijective, hence the same holds  for the composite map:
$$\widetilde Z@>H>> \Cal W\times P^\vee/\Q@>(f^{-1},\,\exp)>>\Cal
I\times Cent(G).$$
\par\noindent{(d)} Restricting  $F$ to $\widetilde Z_{ab}$ induces  a 
surjection
$\widetilde Z_{ab}\twoheadrightarrow \Cal W_{ab}$, and as above bijections
$Z_{ab}\leftrightarrow\Cal W_{ab}$,
$\widetilde Z_{ab}\leftrightarrow \II\times Cent(G)$.
\endproclaim

\vbox{
\demo{Proof}
\par\noindent{(a)}  Assume that $F(z)=F(z')$ for  $z,z'\in\Z$ and set $v=v_z^{-1},\,u=v_{z'}^{-1}$. Then
$v(z+ C_1)=u(z'+
C_1)$, so that $C_1=v^{-1}u(z')-z+v^{-1}u(C_1)$ and
$t_{[v^{-1}u(z')-z]}v^{-1}u\in\Omega$. By Proposition D, (2), we have
either $v=u$ and in turn $z=z'$ (and in this case we are done), or there
exists $j\in J$ such that $v^{-1}u=w_0^jw_0$; moreover
$$z=w_0^jw_0(z')-\omega_j^\vee\tag*$$ which means $z=\psi(z')$ with
$\psi=t_{-\omega_j^\vee}w_0^jw_0\in \Sigma$. Viceversa, if $z=\psi(z'),\,1\ne\psi\in\Sigma$,
relation $(*)$ implies $z+\overline C_1=w_0^jw_0(z'+\overline C_1)$, hence 
$F(z)=F(z')$.
\par
\par\noindent{(b)} In \cite{CP2, Prop.~3} it is proved that
$F_{|Z}:Z\to\Cal W$ is a bijection, with inverse
map
$t_{\tau}v\mapsto v^ {-1}(\tau)$:  by {(a)} we obtain that $F$ maps
the whole $\widetilde Z$ onto $\Cal W$.
\par\noindent{(c)} This statement is a direct consequence of (b) and
Proposition 1.1, (2).
\par\noindent{(d)} This follows from Propositions B, (1), Proposition 1.3, and from 
part (c) of this Proposition.
\endemo
}

\proclaim{Proposition 1.5} If $z\in\widetilde Z$, then
$N(v_z)=\{\a\in \Dp\mid (\a,z)<0\}$.
\endproclaim

\demo{Proof}
We have to prove that the  element $v\in W$ such that
$v(z+C_1)\subset C_\infty$
is defined by the condition $N(v^{-1})=\{\a\in \Dp\mid (\a,z)<0\}$. It suffices
to verify that the element $u$ defined by $N(u^{-1})=\{\a\in \Dp\mid
(\a,z)<0\}$
is such that $(u(z+e),\b)>0\ \forall\,\b\in\Dp,\ \forall\,e\in C_1$.
We can then
conclude that $u=v$ by uniqueness.\par First remark that $N(u)=-uN(u^{-1})$;
then suppose $\b\in N(u)$, or $\b=-u(\gamma),$ $\gamma\in N(u^{-1})$; by
hypothesis $(z,\gamma)\leq -1$, hence we have $(u(z+e),\b)=-(z+e,\gamma)\geq
1-(e,\gamma)>0$. It remains to consider the case $\b\notin N(u)$; in that case
$\b=u(\gamma),\,\gamma\notin N(u^{-1})$ and in particular
$(z,\gamma)\geq 0$, so
that $(u(z+e),\b)=(z+e,\gamma)=(z,\gamma)+(e,\gamma)\geq (e,\gamma)>0$ as
desired.
\endemo
\medskip

\proclaim{Corollary 1.6} If $z\in\Z_{ab}$, then $N(v_z)=\D^{-2}_z\cup
\D^{-1}_z$.
\endproclaim

\proclaim{Proposition 1.7} Suppose $z\in\Z_{ab}$.
Then
$$\i_z=-v_z^{-1}(\D^{-2}_z)\cup v_z^{-1}(\D^{1}_z).$$
\endproclaim

\demo{Proof} Let $w=\widetilde F(z)$, $w=t_\tau v$; thus  
$v_z=v^{-1}$ and $z=v_z(\tau)$.
By Lemma C we obtain that 
$\p_{\i_z}=\D^2_\tau\cup(\D^1_\tau\setminus N(v))$.
Let $\a\in \D^2_\tau$. Since $(v_z(\a), z)=(\a, \tau)=2$ and 
$\D_z^2=\emptyset$, we
obtain $v_z(\a)<0$, thus $-v_z(\a)\in \D^{-2}_z$, or equivalently 
$\a\in -v_z^{-1}(\D^{-2}_z)$.
Conversely, if $\b\in \D^{-2}_z$, then $(\tau, v_z^{-1}(\b))=-2$,
whence $-v_z^{-1}(\b)\in \D^2_\tau$. Therefore 
$\D^2_\tau=-v_z^{-1}(\D^{-2}_z)$.
Similarly, we see that if $\a\in \D^1_\tau\setminus N(v)$, then 
$v_z(\a)\in \D^1_z$, and,
conversely, if $\b\in \D^1_z$, then  $v_z^{-1}(\b)$ belongs to
$\D^1_\tau\setminus N(v)$. This concludes the proof.
\endemo

\remark{\sl Remark}  Corollary 1.6 and Proposition 1.7 show that our
constructions
coincide for points in $\Z_{ab}$ with those performed  by Kostant  in
\cite{Ko2},
taking into account that Kostant's $Z$  is  our $-\Z_{ab}$. In
particular, we have provided
proofs for the results of Sections 2 and 3 and for Proposition 5.2
and Theorem 5.3 of that paper.
\endremark

\medskip
\bigskip
\heading\S2 Compatible Borel subalgebras and compatible discrete series.
\endheading
\bigskip
A symmetric Lie subalgebra of $\g$ is a  subalgebra $\ka$ that is the
fixed point set of an involutory automorphism $\Theta$ of $\g$. If $\ka$ is a
symmetric Lie subalgebra, we set $\lie{p}$ to denote the $-1$ eigenspace of
$\Theta$. It follows that we have a decomposition
$$
\g=\ka\oplus\lie{p}
$$
that is usually referred to as a Cartan decomposition of $\g$.
     An equal rank sym\-me\-tric Lie
subalgebra is a symmetric Lie subalgebra that contains a Cartan subalgebra of
$\g$. A procedure to construct equal rank symmetric Lie subalgebras is the
following:  if $\tau\in P^\vee$,
set $\Theta_\tau=Ad(\exp(\sqrt{-1}\pi\tau))$. Then $\Theta_\tau$ is
an involutory
automorphism of $\g$. If we set
$$
\g=\ka_\tau\oplus \lie{p}_\tau
$$
to be the corresponding Cartan decomposition, then $\ka_\tau$ is an equal rank
symmetric Lie subalgebra.

On the other hand, if $\ka$ is an equal rank symmetric Lie subalgebra of $\g$
and  $\h'$ is a Cartan subalgebra contained in $\ka$, then, by \cite{Hel},
Ch.~IX, Proposition~5.3, there is
$x\in
\g$ such that
$\Theta=Ad(\exp(x))$. Moreover, by Theorem~5.15 of Ch.~X of \cite{Hel}, there
is an automorphism $\phi$ of $\g$ such that $\phi\Theta\phi^{-1}$ is an
automorphism $\Theta'$ of type $(s_0,\dots,s_r;k)$ for $\h$, $r$ being
    the rank of the subalgebra of $\Theta'$-fixed points. Since $\ka$ is
an equal rank symmetric subalgebra, we have $r=n$. Since
$\Theta'=Ad(\exp(\phi(x)))$, it follows from Theorem~5.16 (i),
Ch.~X of \cite{Hel} that  $k=1$, i.\ e.  $\Theta'$ fixes $\h$
pointwise.

Set $X=\overline{C}_2\cap P^\vee$.
If we set $\tau=\sum_i s_i\omega_i^\vee$, 
then $\tau\in X$. Indeed, since $\Theta'$ is an involution,
then, according
to Theorem~5.15 of \cite{Hel} again,
$$
2=s_0+\sum_{i=1}^n m_is_i
$$
(recall that $m_i$ is the coefficient of $\al_i$ in the highest root).

Note that $\Theta'=\Theta_\tau$. In fact, calculating the action on
the root vectors $X_{\a_i},\,1\leq i\leq n,\ X_{-\th}$
(which are Lie algebra generators for $\g$), we have:
$$\align
&\Theta_\tau(X_{\al_i})=(-1)^{s_i}X_{\al_i}=\Theta'(X_{\al_i})\\
&\Theta_\tau(X_{-\theta})=(-1)^{s_0-2}X_{-\theta}=(-1)^{s_0}X_{-\theta}=
\Theta'(X_{-\theta}).\endalign$$

Finally we observe that, by Theorem~5.4 of Ch. IX of \cite{Hel}, we can write
$\phi=\nu Ad(g)$, with $\nu$ an automorphism of $\g$ leaving $\h$ and
$C_\infty$ invariant. Then
$$
\nu^{-1}\Theta_\tau\nu=\Theta_{\nu^{-1}(\tau)}=Ad(g)\Theta Ad(g^{-1}).
$$
Since $\nu^{-1}(\tau)\in X$, we have proved a weaker form of Proposition~4.1 of
\cite{Ko2}, namely:
\proclaim{Theorem 2.1}
If $\ka$ is an equal rank symmetric Lie subalgebra of $\g$ then there is
$\tau\in X$ such that $\ka_\tau=Ad(g)\ka$ for some $g\in G$.
\endproclaim

We now turn our attention to the $\Theta$-stable Borel subalgebras. Suppose
that $\ka$ is an equal rank symmetric subalgebra of $\g$ and let $\Theta$ be
the corresponding involution.

We let $K_\tau$
be the subgroup of $G$ corresponding to $\ka_\tau$.

Theorem~1 of \cite{Ma} gives the following characterization of
$\Theta$-stable Borel subalgebras.

\proclaim{Theorem 2.2}
     If
$\lie{b}'$ is a $\Theta$-stable Borel subalgebra of $\g$, then there exist
$g\in G$,  $w\in W$, and  $\tau\in X$ such
that
$Ad(g)\ka=\ka_\tau$ and
$Ad(g)\lie{b}'=w\lie{b}$.
\endproclaim

\demo{Proof}
By Theorem~2.1, we can find $g'\in G$ such that
    $$
Ad(g')\ka=\ka_\tau
\hbox{ and }
Ad(g')\Theta Ad(g'{}^{-1})=\Theta_\tau,
$$
so $Ad(g')\lie{b}'$ is
$\Theta_\tau$-stable. By Theorem~1 of
\cite{Ma}, we can find an element $k\in K_\tau$ such that
$\lie{b}''=Ad(kg')\lie{ b}'$ contains a $\Theta_\tau$-stable Cartan subalgebra
$\h'$.
Set $\Delta'$ to be the set of roots of $(\g,\h')$ and let $\Delta'{}^+$
denote the positive system in $\Delta'$ defined by $\lie{b}''$.
Since
$\lie{b}''$ and $\h'$ are
$\Theta_\tau$-stable, it follows that the map $\al\mapsto\al\circ \Theta_\tau$
defines an automorphism of the Dynkin diagram. Since $\Theta_\tau$ is of inner
type we conclude that $\Theta_\tau$ fixes pointwise $\h'$, i. e.
$\h'\subset\ka_\tau$. Hence there is
$k'\in K_\tau$ such that $Ad(k')\h'=\h$.
Set $g=k'kg'$; then $Ad(g)\ka=\ka_\tau$
and $Ad(g)\lie{b}'$ is a Borel subalgebra containing $\h$. Then $Ad(g)\lie{b}'$
defines a positive system in $\Delta$, therefore there is an element $w\in W$
such that $Ad(g)\lie{b}'=w\lie{b}$.
\endemo

\proclaim{Definition} Let $\g=\ka\oplus\lie{p}$ be a Cartan
decomposition, with Cartan involution $\Theta$. We will say that a
Borel subalgebra
$\lie{b}'$ is compatible with
$\Theta$ (or with
$\ka$) if it is $\Theta$-stable and
$$
[[\lie{b}'_{\lie{p}},\lie{b}'_{\lie{p}}],\lie{b}'_{\lie{p}}]=0,
$$
where $\lie{b}'_{\lie{p}}=\lie{b}'\cap{\lie{p}}$.\endproclaim
Clearly $\lie{b}$ is compatible with $\ka_\tau$ for any $\tau\in X$.
Conversely we have the following theorem.
\proclaim{Theorem 2.3} \cite{Ko2, Theorem~4.3}
Let $\ka$ be an equal rank symmetric subalgebra such that
$\lie{h}\subset\lie{k}$. If
$\lie{k}$ is compatible with
$\lie{b}$ then $\ka=\ka_\tau$ for some $\tau\in X$. Moreover if $\ka'$ is
any equal rank symmetric subalgebra, then a Borel subalgebra
$\lie{b}'$ is compatible with
$\ka'$ 		if and only if there exist
$g\in G$ and
$\tau\in
X$ such that $Ad(g)\ka'=\ka_{\tau}$ and
$Ad(g)\lie{b}'=\lie{b}$.
\endproclaim
\demo{Proof}
Suppose that $\lie{k}$ is compatible with $\lie{b}$. Since $\lie{h}\subset\ka$,
by Exercise C.3 of Ch.~IX of \cite{Hel}, we have that $\Theta=Ad(\sqrt{-1}\pi
h)$ with $h\in P^\vee$. Clearly we can choose $h=\sum\epsilon_i\omega_i^\vee$
with $\epsilon_i\in \{0,1\}$.  We claim
that, since
$\lie{b}$ is compatible, then
$\sum m_i\epsilon_i\le 2$, i.e. $h\in X$.

Indeed write
$$
\bb_0=\h\oplus\sum_{\al(h)=0}\g_\al,\quad\bb_1=\sum_{\al(h)=1}\g_\al.
$$
Clearly
$$
[\bb_0+\bb_1,\bb_0+\bb_1]\subset \bb_0+\bb_1+[\bb_1,\bb_1]
$$
hence, since $\bb$ is compatible,
$$
[[\bb_0+\bb_1,\bb_0+\bb_1],\bb_0+\bb_1]\subset \bb_0+\bb_1+[\bb_1,\bb_1].
$$
Since $\bb_0+\bb_1$ generates $\bb$ we have that
$
\bb=\bb_0+\bb_1+[\bb_1,\bb_1]$, thus, in particular, the root vector
$X_{\theta}$ belongs to either $\bb_0$, or $\bb_1$, or $[\bb_1,\bb_1]$.

If $X_{\theta}\in \bb_0$ then $\theta(h)=0$, if $X_{\theta}\in
\bb_1$ then $\theta(h)=1$, and, if $X_{\theta}\in [\bb_1,\bb_1]$, then
$\theta(h)=2$. Since $\theta(h)=\sum\epsilon_im_i$ the first result
follows.

For the second assertion we use Theorem~2.2 to deduce that there is
$g\in G$ such that $Ad(g)\bb'=\bb$ and $Ad(g)\ka'=\ka_{w\tau'}$ for some
$\tau'\in X$ and $w\in W$. By the first part of the proof we obtain
$\ka_{w\tau'}=\ka_{\tau}$ for some $\tau\in X$.
\endemo

Given $z\in \h_\R$, let $v_z$ be the unique element of $W$ such that
$$N(v_z)=\{\al\in\Delta^+\mid (\al,z)<0\}.$$
 This definition coincides for $z\in \widetilde Z$ with the one already
introduced (see Proposition~1.5). \par
It is easy to check that, if $z\in \h_{\R}$, then $v_z^{-1}(z)$ is dominant,
hence we can define
$$dom:\h_\R\to \overline{C}_\infty \quad z\mapsto v^{-1}_z(z).$$
If $z\in \widetilde{Z}$ and $\widetilde F(z)=t_\tau v$, then $dom(z)=\tau$.
Clearly
$dom(\widetilde{Z}_{ab})= X$, in fact $w_{0}\tau\in\Z_{ab}$
for all $\tau\in X$.
\par
If $\tau\in X$, set
$\widetilde{Z}_\tau=dom^{-1}(\tau)\cap \widetilde{Z}_{ab}$.
Let $W_\tau$ be the Weyl group of $(\ka_\tau,\h)$ and
denote by
$n_\tau$  the index of $W_\tau$ in $W$.
The following result affords a proof of  \cite{Ko2, Theorem~4.5}.
\proclaim{Theorem 2.4} The map $z\mapsto W_{\tau}v_{z}^{-1}$ is a bijection between
$\widetilde{Z}_{\tau}$ and $W_{\tau}\backslash W$. In particular, 
for any $\tau\in X$, one has
$$
|\widetilde{Z}_\tau|=n_\tau.
$$
\endproclaim
\demo{Proof}
We note that, if  $z\in \widetilde{Z}_{\tau}$ then
$$
\Delta^2_\tau=-v_z^{-1}(\Delta^{-2}_z),
$$
$$
\Delta_\tau^0=v_z^{-1}(\Delta^0_z),
$$
$$
\Delta_\tau^1=v_z^{-1}(\Delta^1_z)\cup -v_z^{-1}(\Delta_z^{-1}).
$$

The root system for
$(\ka_\t,\h)$ is $\Delta_\tau=\pm\Delta_\tau^0\cup\pm \Delta^2_\tau$. By
our formulas above
$ v_z^{-1}(\Delta^+)$ contains $-\Delta^2_\tau\cup \Delta^{0}_{\tau}$.
This implies that  $W_\tau v_z^{-1}=W_\tau v_{z'}^{-1},\ z,z'\in \widetilde Z_\t,$ if and only if $z=z'$:
indeed if $v_z^{-1}=w'v_{z'}^{-1}$ for some $w'\in
W_\tau$ then
$$
-\Delta^2_\tau\cup \Delta^{0}_{\tau}\subset v_z^{-1}(\Delta^+)
=w'v_{z'}^{-1}(\Delta^+)\supset w'(-\Delta^2_\tau\cup \Delta^{0}_{\tau}).
$$
Since $-\Delta^2_\tau\cup \Delta^{0}_{\tau}$ is a positive system for
$\Delta_{\tau}$, it follows that $w'=1$. Hence $v_z=v_{z'}$ and in turn
$z=v_z(\tau)=v_{z'}(\tau)=z'$.

We now verify that  $W_\tau w=W_\tau v_z^{-1}$ for some
$z\in\widetilde{ Z}_\tau$:
let $w'\in W_{\tau}$ be the unique element such that
$w'(w(\Delta^{+})\cap\Delta_{\tau})=-\Delta^{2}_{\tau}\cup
\Delta_{\tau}^{0}$ so that we have that
$w'w(\Delta^{+})\supset-\Delta^{2}_{\tau}\cup\Delta^{0}_{\tau}$.

Set $z=(w'w)^{-1}(\tau)$.
Let us verify that $z\in \widetilde{Z}_{ab}$: if $\al\in \Delta^{+}$
then $\al(z)=w'w\al(\tau)\in
\{\pm2,\pm1,0\}$, so it is enough to verify that
$\Delta^{2}_{\tau}$ is disjoint from $w'w(\Delta^{+})$, but this is obvious
since
$w'w(\Delta^{+})$
contains $-\Delta_{\tau}^{2}$.

We are left with showing that $W_{\tau}w=W_{\tau}v_{z}^{-1}$: since $z$
is in the $W$-orbit of $\tau$, then $v_{z}^{-1}(z)=\tau=w'w(z)$,
hence  $w'wv_z\in Stab_W(\t)\subset W_\t$; therefore $W_\t w=W_\t v_z^{-1}$.
\endemo
We are going to prove a result which implies \cite{Ko2, Theorem~5.5}.
As shown in the
proof of Theorem~2.4,
$\Delta_\tau=\pm\Delta_\tau^0\cup\pm\Delta^2_\tau$ is the set of roots for
$\ka_\tau$ and  $\Delta_\tau^+=-\Delta^2_\tau\cup\Delta^0_\tau$ is a positive
system for $\D_\tau$. We let $\Pi_\tau$ denote the set of simple roots for
$\D_\tau$ corresponding to $\D^+_\tau$.
\smallskip
\def\stab{{\text{Stab}}}

We consider $\widehat W_2=W \ltimes T(2 Q^\vee)$: $\widehat W_2$
can be viewed as the affine Weyl group of ${1\over 2} \D$, so
that $C_2$ is its fundamental alcove relative to ${1\over 2} \Pi$.
In the following lemma we describe $\Pi_\tau$  and see that $W_\tau$
is strictly related with
the stabilizer $\stab_{\widehat W_2}(\tau)$ of $\tau$ in $\widehat W_2$.

\proclaim{Lemma~2.5}
(1). We have $\stab_{\widehat W_2}(\tau)\subset W_\tau\ltimes T(2 Q^\vee)$.
Moreover, $\stab_{
\widehat W_{2}}(\tau)\cong W_\tau$ via the canonical
projection
of $\widehat W$ onto $W$.
\newline
(2).
Let $\Pi_\tau$ denote the set of simple roots for $\D_\tau$ relative to
$\D^+_\tau$. If $(\tau, \theta)< 2$, then $\Pi_ \tau=\Pi\cap \tau^\bot$;
if $(\tau, \theta)=2$, then $\Pi_\tau=(\Pi\cap \tau^\bot)\cup\{-\theta\}$.
\endproclaim

\demo{Proof}
We start with some general facts. Consider any point $\nu\in\overline C_1$.
It is well known (see \cite{Bou, V, 3.3}) that the stabilizer
$\stab_{\widehat
W}(\nu)$ of $\nu$ in $\widehat W$ is the parabolic subgroup generated by
the
simple reflections which fix $\nu$.
Set $\D(\nu)=\{\a\in \D\mid (\a, \nu)\in \ganz\}$,
$\widehat \D_\nu= \{\a-(\a, \nu)\d\mid \a\in \D(\nu)\}$,
and $\widehat \Pi_\nu=\widehat \D_\nu\cap \widehat \Pi$.
It is clear that both $\D(\nu)$ and $\widehat \D_\nu$ are root subsystems
of $\D$ and $\widehat \D$, respectively; moreover, $\widehat \D_\nu$
is isomorphic to $\D(\nu)$, via the natural projection on $\frak
h^*_\real$.
Since $\nu\in\overline C_1$, we have that $(\a, \nu)\in \{0,-1, 1\}$ for
each
$\a\in \D(\nu)$, hence we easily obtain that $\widehat \D_\nu$ is the
standard
parabolic subsystem generated by $\widehat \Pi_\nu$.
Since the simple reflections that fix $\nu$ are exactly the reflections
corresponding to $\widehat \Pi_\nu$, we have that $\stab_{\widehat
W}(\nu)$ is
the Weyl group of $\widehat \D_\nu$.
This implies in particular that the canonical projection of
$\stab_{\widehat W}(\nu)$ on $W$ provides an isomorphism of
$\stab_{\widehat W}(\nu)$ with the Weyl group of $\D(\nu)$.
Moreover, if we project $\widehat \Pi_\nu$ on $\D$, we obtain a basis
for $\D(\nu)$. Thus we have that: if $(\nu, \theta)<1$, then
$\Pi\cap \nu^\bot$ is a basis of $\D(\nu)$; if $(\nu, \theta)=1$,
then $(\Pi\cap \nu^\bot)\cup \{-\theta\}$ is a basis of $\D(\nu)$.
\par
If we apply the above remarks to $\widehat W_2$ in place of
$\widehat W$, $C_2$ in place of $C_1$,  and $\tau \in X$ in place of
$\nu$, we obtain
that
$\stab_{\widehat W_2}(\tau)$ is a parabolic subgroup of $\widehat W_2$
which projects isomorphically onto
$W_\tau$ via the canonical projection of $\widehat W_2$ on $W$.
In particular, we have
$\stab_{\widehat W_2}(\tau)\subset W_\tau \ltimes T(2Q^\vee)$.
Moreover we have that if $(\tau, \theta)< 2$, then
$\Pi\cap \tau^\bot$ is a basis of $\D_\tau$; if  $(\tau, \theta)=2$, then
a basis
of $\D_\tau$ is given by  $\{-\theta\} \cup (\Pi\cap \tau^\bot)$.
\qed\enddemo

\smallskip
  Recall  that
$\Omega_2=\left\{t_{2\omega_j^\vee}w_0^jw_0\mid j\in
J\right\}\cup\{1\}$ and define, for $z\in \widetilde Z$, 
$$\lie{b}_z=v_z^{-1}\lie{b}.$$

\proclaim{Theorem 2.6} Suppose that
$\tau\not\in\{2\omega^\vee_i\mid i\in J\}\cup\{0\}$.
Set
$$
\widetilde{Z}_\tau^{\text{cmpt}}=\{z\in\widetilde{Z}_\tau\mid
\lie{b}_z \hbox{ is
compatible with\ }\ka_\tau\}.
$$
Then there is a canonical bijection between
$\widetilde{Z}_\tau^{\text{cmpt}}$ and
$\Omega_2$. In particular
$$
|Cent(G)|=|\widetilde{Z}_\tau^{\text{cmpt}}|.
$$
\endproclaim
\demo{Proof}
If $s\in\Omega_2$ write $s=t_\nu v$ with $v\in W$ and $\nu\in 2P^\vee$. By
Theorem~2.4 there is a unique $z_s\in\widetilde{Z}_\tau$ such that
$W_\tau v_{z_s}^{-1}=W_\tau v$. Notice that $z_s\in\Comp$: indeed, if
$\tau'=s^{-1}(\tau)$, we have that $\t'\in X$, hence   $\bb$ is compatible with
$\ka_{\tau'}=\ka_{s^{-1}(\tau)}=\ka_{v^{-1}(\tau)}$. In turn  $v\bb$ is compatible
with $\ka_\tau$, but $v_{z_s}^{-1}=w'v$ with $w'\in W_\tau$, hence
$\bb_{z_s}=w'v\bb$ is compatible with $w'\ka_\tau=\ka_\tau$.
We can thus define a map $B:\Omega_2\to\Comp$ by setting $B(s)=z_s$.

Let us show that $B$ is injective: we need to show that, if
$s=t_\nu v$ and
$s'=t_{\nu'}v'$ are such that $W_\tau v=W_\tau v'$, then $v=v'$. Since
$\Omega_2$ is a group, it is enough to check that, if $s=t_\nu v$ is such that
$v\in W_\tau$, then $v=1$. The assumption that $v\in W_\tau$ implies that
$v(\Delta_\tau)=\D_\tau$ and, by \cite{IM, Prop.~1.26, (ii)},
$v(\Pi\cup\{-\theta\})=\Pi\cup\{-\theta\}$.
Since $\tau\not\in\{2\omega^\vee_i\mid i\in J\}\cup\{0\}$, we see that
$\Pi_\tau=\Delta_\tau\cap(\Pi\cup\{-\theta\})$, hence $v(\Pi_\tau)=\Pi_\tau$
and $v=1$.

It remains to show that $B$ is surjective. Fix $z\in\Comp$: since
$v_z^{-1}\bb$ is compatible with $\ka_\tau$ it follows that $\bb$ is compatible
with $\ka_{v_z(\tau)}=\ka_z$. By Theorem~2.3, we deduce that there is $\tau'\in
X$ such that $\ka_z=\ka_{\tau'}$, hence $z=\tau'+\nu$ with $\nu\in 2P^\vee$,
or, equivalently,
$$
\tau=v_z^{-1}(\tau'+\nu)=t_{\nu'}v_z^{-1}(\tau),
$$
where $\nu'=v_z^{-1}(\nu)$.

Since $\overline{C}_2$ is a fundamental
domain  for $\widehat W_2$, there is a unique element $\tilde{u}\in \widehat
W_2$ such that $\tilde{u}t_{\nu'}v_z^{-1}(\overline{C}_2)=\overline C_2$. Set
$s=\tilde{u}t_{\nu'}v_z^{-1}$: clearly $s\in\Omega_2$. Since $\tau\in
t_{\nu'}v_z^{-1}(\overline{C}_2)\cap\overline C_2$, we find that
$\tilde{u}\tau=\tau$, therefore, by the previous Lemma
    we have
$\tilde{u}\in W_\tau\ltimes T(2\Q)$. Thus $B(s)=z$ and we are done.
\endemo

Fix again $\tau\in X$ and set $\sigma$ to denote a conjugation in
$\g$ corresponding to a compact real
form $\g_u$ such that $\Theta_\tau\sigma=\sigma\Theta_\tau$ and set
$\sigma_\tau=\Theta_\tau\sigma$.
Then
$\sigma_\tau$ is a conjugation of $\g$ defining a real form $\g_\tau$. We set
$G_\tau$ to be the subgroup of
$G$ corresponding to $\g_\tau$ and $K$ to be the subgroup of $G_\tau$
corresponding to $\g_\tau\cap\ka_\tau$.

Set $P_{\text{reg}}=\{\lambda\in P\mid (\lambda,\al)\ne 0\hbox{ for
all } \al\in\Delta\}$. Let $\widehat G_{\tau}^{\text{disc}}\subset\widehat
G_{\tau}$ denote the set of equivalence classes of discrete series for
$G_{\tau}$. If $\lambda\in P_{\text{reg}}$  we let $\pi_\lambda\in
\widehat G_{\tau}^{\text{disc}}$
denote the equivalence class
corresponding to the parameter $\lambda$ in Harish
Chandra parametrization. We recall that given $\lambda,\mu\in
P_{\text{reg}}$ then
$\pi_\lambda=\pi_\mu$  if and only if there is
$w\in W_\tau$
such that $w\lambda=\mu$.

If $\lambda\in P_{\text{reg}}$ then we set
$\D^+_\lambda=\{\al\in\D\mid (\lambda,\al)>0\}$ and
$\bb_\lambda=\h\oplus\sum_{\al\in\D^+_\lambda}\g_\al$ to be the
corresponding Borel subalgebra.
Given
$z\in
\widetilde{Z}_\tau$, we set $P_{\text{reg}}(z)=\{\lambda\in P_{\text{reg}}\mid
\bb_\lambda=\bb_z\}$ and  $\rho_z=v_z^{-1}\rho$. Set also $\widehat
G^{\text{disc}}_\tau(z)=\{\pi_\lambda\mid \lambda\in
P_{\text{reg}}(z)\}$. As observed in \cite{Ko2}, Theorem~2.4 implies
that the sets $\widehat
G^{\text{disc}}_\tau(z)$ give a partition of the set of
equivalence classes of discrete
series for $G_\tau$, namely
$$\widehat G_\tau^{\text{disc}}=\bigcup_{z\in Z_\tau}\widehat
G_\tau^{\text{disc}}(z).$$ 

With the notation of \S 1, if
$z\in\widetilde{Z}_{ab}$, recall that 
$\i_z\in\I_{ab}$  denotes
$f^{-1}(F(z))$ and
$\Phi_z\equiv\Phi_{\i_z}=\{\alpha\in\Delta^+\mid \g_\a\subset \i_z\}$. Set also
$\Phi^2_z=-v_z^{-1}(\Delta^{-2}_z)$, $\Phi^1_z=v_z^{-1}(\Delta^{1}_z)$. By
Proposition~1.7 we know that $\Phi_z=\Phi^2_z\cup \Phi_z^1$.
\smallskip

The relation between abelian ideals and discrete series is not as good as one
might expect, for example the correspondence depends on the choice of $\tau$
even when different choices for $\tau$ give rise to isomorphic real forms of
$G$. Nevertheless Kostant observed in
\cite{Ko2} that there are  relations between the elements of
$\widetilde{Z}_\tau$ and
the structure of the corresponding  discrete series. The first occurrence of
such a relation involves the realization of the discrete series via
$L^2$-cohomology: one can give to the principal bundle
$G_\tau/T$   a complex structure by declaring that
$\lie n$ is the space of antiholomorphic tangent vectors; if
$\lambda\in P_{\text{reg}}(z)$ then set $\Cal L_{\lambda-\rho}$ to be the
holomorphic line bundle on $G_\tau/T$ associated to the character
$\exp(\lambda-\rho)$ of $T$. According to \cite{S2}, the
$L^2$-cohomology groups
$H^p(\Cal L_{\lambda-\rho})$ vanish except  in one degree $k(\lambda)$ and
$H^{k(\lambda)}(\Cal L_{\lambda-\rho})$ gives a realization of $\pi_\lambda$.
It turns out that the decomposition $\Phi_z=\Phi^2_z\cup \Phi_z^1$
implies, applying Theorem~1.5 of \cite{S2}, that, if
$\lambda\in P_{\text{reg}}(z)$, then
$k(\lambda)=\text{dim}(\i_z)$ (cf. \cite{Ko2, Theorem 6.5}).

Another instance, discussed in \cite{Ko2}, where there is a connection between
abelian ideals and discrete series involves the $K$-structure of $\pi_\lambda$.
   Recall that, if $\lambda\in
P_{\text{reg}}(z)$, then
     $\Delta_\tau^+=-\Delta^2_\tau\cup\Delta^0_\tau$ is a positive
system for $\Delta_\tau$ contained in $\Delta^+_\lambda=v_z^{-1}\Delta^+$.
Set $\rho_c=\half\sum_{\al\in\D^+_\tau}\al$ and $\rho_n=\rho_z-\rho_c$.
It is well known that the minimal $K$-type of $\pi_\lambda$ has
highest weight with respect
to $\D^+_\tau$ given by the formula
$$
\mu_\lambda=\lambda+\rho_n-\rho_c.
$$
We are going to prove, using this formula, that if    
 $\lambda\in P_{\text{reg}}(z)$ then the highest weight
$\mu_\lambda$ of the minimal
$K$-type of $\pi_\lambda$ with respect to $\Delta^+_\tau$ is
$$
\mu_\lambda=\lambda-\rho_z+2\langle\i_z\rangle-\half\tilde{\tau}.
$$
Notation is as follows: for  $h\in\h_\R$ we set
$
\tilde{h}=\sum_{\alpha\in\Delta}(\al,h)\al.
$
and for $\i\in\I$ we put $\langle \i\rangle=\sum\limits_{\a\in \p_i}\a$.
A  special  case of the previous formula is
$$
\mu_{\rho_z}=2\langle\i_z\rangle-\half\tilde{\tau}.
$$
To prove the  formula, just compute:
$$
\align
\mu_\lambda&=\lambda+\rho_n-\rho_c\\
&=\lambda-\rho_z+2\rho_n\\
&=\lambda-\rho_z+\sum_{\al\in\Delta^1_z}v_z^{-1}\al+
\sum_{\al\in\Delta^{-1}_z}v_z^{-1}\al\\
&=\lambda-\rho_z+\sum_{\al\in\Phi^1_z}\al-
\sum_{\al\in\Delta^1_\tau\backslash\Phi^1_z}\al\\
&=\lambda-\rho_z+2\sum_{\al\in\Phi_z}\al-\sum_{\al\in\Phi^1_z}\al
-\sum_{\al\in\Delta^1_\tau\backslash\Phi^1_z}\al-2\sum_{\al\in\Phi^2_z}\al\\
&=\lambda-\rho_z+2\langle\i_z\rangle-\sum_{\al\in\Delta^1_\tau}\al
-\sum_{\al\in\Delta^2_\t}2\al\\
&=\lambda-\rho_z+2\langle\i_z\rangle-\sum_{\al\in\Delta^+}\al(\tau)\al\\
&=\lambda-\rho_z+2\langle\i_z\rangle-\half\tilde{\tau}.
\endalign
$$
The previous calculation proves Theorem 6.6 of \cite{Ko2}.\smallskip

We now turn our attention to the discrete series that correspond to elements of
$\widetilde{Z}_\tau^{\text{cmpt}}$.

\proclaim{Definition}
A discrete series  $\pi_{\lambda}$ is said compatible if $\lambda\in
P_{\text{reg}}(z)$ with
$z\in\widetilde{Z}^{{\text{cmpt}}}_{\tau}$.\endproclaim
By
Theorem~2.6,
the number of compatible discrete series is $|P^{\vee}/Q^{\vee}|$.
In particular it is independent
from the particular real form $G_{\tau}$.

Following \S~8 of Ch. XI of \cite{KV},
we recall how a representative of $\pi_\lambda$ is constructed: set
$$
\lie{n}_{\lambda}=[\bb_\lambda,\bb_\lambda],
$$
$S=\dim \lie{n}_\lambda\cap\ka$, and
$\rho(\lie{n}_\lambda)=\half\sum_{\alpha\in\Delta_\lambda^+}\al$. If
$\mu\in P$ we denote by
$\C_\mu$  the one dimensional representation of $\bb_\lambda$ defined
by setting
$(h+n)\cdot c=\lambda(h)c$.

    Then, according
to Theorem~11.178 of \cite{KV}, the $(\g,K)$-module of a representative of
$\pi_\lambda$ is
$$
V^\lambda_{(\g,K)}=
({}^u\DR^{(\g,K)}_{(\bb_\lambda,T)})^S
(\C_{\lambda+\rho(\lie{n}_\lambda)}).
$$
Recall that
$({}^u\DR^{(\g,K)}_{(\bb_\lambda,T)})^j=(\Gamma_{(\g,T)}^{(\g,K)})^j
\circ \text{pro}_{(\bb_\lambda,T)}^{(\g,T)}$, where
$(\Gamma_{(\g,T)}^{(\g,K)})^j$ is the $j$-th derived functor of the Zuckerman
functor, while $\text{pro}_{(\bb_\lambda,T)}^{(\g,T)}$ is the ordinary
algebraic induction functor:
$\text{pro}_{(\bb_\lambda,T)}^{(\g,T)}(Z)=\text{Hom}_{\bb_\lambda}(U(\g),Z)_{T}$.

     In what follows we adopt the following notation: if $\lie{a}$ is any
subspace of $\g$ that is stable for the action of $\h$, we let
$\Delta(\lie{a})$ be the set of roots occurring in $\lie{a}$, in other words
$\al\in\Delta(\lie{a})$ if and only if $\g_\al\subset \lie{a}$.

     Set
$\Pi_\lambda$ to be the set of simple roots for $\Delta_\lambda^+$ and
$\Pi_\lambda^c\subset\Pi_\lambda$ to be the set of simple compact roots. Let
$$
\lie{ q}=\lie{ m}\oplus\lie{ u}\eqno(1)
$$
be the
corresponding parabolic; i.e. $\Delta(\lie{ m})$ is the
subsystem of $\Delta$ generated by $\Pi_\lambda^c$ and
$\Delta(\lie{ u})=\Delta_\lambda^+\,\backslash\,(\Delta(\lie{
m})\cap \Delta_\lambda^+)$.

Set $\bb_{\lie{m}}=\bb_\lambda\cap\lie{m}$, $\D^+(\lie m)=\Delta(\lie{ m})\cap
\Delta^+_{{\lambda}}$,
$\rho_{\lie
{m}}=\half\sum_{\alpha\in \Delta^+(\lie{ m})}\alpha$, $M$  the subgroup of $K$
corresponding to $\lie{m}\cap\g_\tau$, and
$$
\lambda_1=\lambda+\rho(\lie n_\lambda)-\rho_{\lie{m}}.\eqno(2)
$$
Then $\lambda_1$ is dominant and regular for
$\Delta^+(\lie{m})$ and, according to
Corollary~4.160 of \cite{KV},
$$
({}^u\DR^{(\lie{m},M)}_{(\bb_{\lie{m}},T)})^q
(\C_{\lambda+\rho(\lie{n}_\lambda)})=\cases
0&\hbox{if $q\ne\dim \bb_{\lie{m}}\cap\lie u$}\\
V_{(\lie m,M)}^{\lambda_1}&\hbox{if $q=\dim \bb_{\lie{m}}\cap\lie u$}
\endcases\eqno(3)
$$

Let  $s=\dim (\lie{u}\cap\ka)$. Combining
(3) with Corollary~11.86 of \cite{KV}, we find that
$$
V^\lambda_{(\g,K)}=({}^u\DR^{(\lie{g},K)}_{(\lie{q},M)})^s
(V_{(\lie m,M)}^{\lambda_1}).\eqno(4)
$$

Fix $z\in\widetilde{Z}_\tau^{\text{cmpt}}$. We are going to compute the
$K$-spectrum of the compatible discrete series, that is the restriction
to $K$ of $V^{\lambda}_{(\g,K)}$ when
$\lambda\in P_{\text{reg}}(z)$.
Set
$S(\lie u\cap\lie p)$ to denote the symmetric algebra of $\lie u\cap\lie p$.

\proclaim{Lemma 2.7}Suppose that $z\in\Comp$ and $\lambda\in
P_{\text{reg}}(z)$.
If  $\mu$ is
$\Delta^+(\lie{m})$-dominant and  there is $n\in\nat$ such that
$$\dim\hbox{\rm Hom}_M(V_{(\lie m,M)}^{\mu+\rho_{\lie m}}, S^n(\lie u\cap\lie
p)\otimes V_{(\lie m,M)}^{\mu_\lambda+\rho_{\lie m}})\ne 0$$ then
$\mu$ is dominant for
$\Delta_\tau^+$.
\endproclaim
\demo{Proof}
Since
$\lie{b}_z$ is compatible with $\ka_\t$, then $\lie{b}$ is compatible with
$\ka_z$, hence, by Theorem~2.3, $\ka_z=\ka_{\tau'}$ for some
$\tau'\in X$.
As in the proof of Theorem~2.3, we can write
$$
\bb=\bb_0+\bb_1+[\bb_1,\bb_1],
$$
where $\bb_0=\h\oplus\sum_{\al\in \D^+\atop \al(\tau')=0}\g_\al$ and
$\bb_1=\sum_{\al(\tau')=1}\g_\al$.
It is then clear that $\lie m\cap\bb_z\supset v_z^{-1}\bb_0$ and
$v_z^{-1}\bb_1=\lie u\cap
\lie p$. This implies that $v_z^{-1}[\bb_1,\bb_1]\subset\lie u\cap\ka$. Since
$\bb_z=v_z^{-1}\bb_0+v_z^{-1}\bb_1+v_z^{-1}[\bb_1,\bb_1]$, we deduce that $\lie
u\cap\ka=v_z^{-1}[\bb_1,\bb_1]$ and $\lie m\cap\bb_z=v_z^{-1}\bb_0$.
Hence $[\lie u\cap\lie
p,\lie u\cap\ka]=0$ and $\lie u\cap\ka$ is abelian. In particular, we have that
$(\al,\beta)\ge 0$ whenever $X_\al\in\lie u\cap\lie p$ and $X_\beta\in\lie
u\cap\ka$.

Since  $\mu_\lambda$ is dominant for $\Delta^+_\tau$,
    our result follows easily as in Lemma~3.1 of \cite{EPWW}.
\endemo

\remark{\sl Remark}
What we actually observe in Lemma~2.7 is the fact that the compatible
discrete series of
\cite{Ko2} are exactly the small discrete series of \cite{GW}. In
\cite{EPWW} it is also
shown that Lemma~2.7 implies that one can compute the full
$K$-spectrum of the discrete
series. In the next result we prove this fact using directly
Blattner's formula.
\endremark

\proclaim{Theorem 2.8}
If $z\in \Comp$, $\l\in P_{\text{reg}}(z)$ and $\mu$ is a
$\Delta_\tau^+$-dominant weight, then
$$
\dim\hbox{\rm Hom}_{K}(V_{(\ka,K)}^{\mu+\rho_c}
,V_{(\g,K)}^{\lambda})=\sum_{n=0}^{+\infty}
\dim\hbox{\rm Hom}_{M}(V_{(\lie m,M)}^{\mu+\rho_{\lie m}},S^n(\lie u\cap\lie
p)\otimes V_{(\lie m,M)}^{\mu_\lambda+\rho_{\lie m}}).\eqno(5)
$$
\endproclaim
\demo{Proof}
Write
$$
m_\mu=\dim\hbox{\rm Hom}_{K}(V_{(\ka,K)}^{\mu+\rho_c}
,V_{(\g,K)}^{\lambda}).
$$

    As in
(11.73) of \cite{KV}, we write
$$
({}^u\DR^{(\lie{g},K)}_{(\lie{q},M)})^j
(Z)=\DR^j(Z\otimes(\Lambda^{\text{top}}\lie u)^*),
$$
so, by (4),
$$
\DR^q
(V_{(\lie m,M)}^{\lambda_1}\otimes(\Lambda^{\text{top}}\lie u)^*)=\cases
0&\hbox{if $q\ne s$}\\
V_{(\lie g,K)}^{\lambda}&\hbox{if $q=s$}
\endcases.
$$

Applying Theorem 5.64 of \cite{KV} we find that
$$m_\mu=\sum_{j=0}^s(-1)^{s-j}\sum_{n=0}^\infty\dim\hbox{\rm Hom}_M(H_j(\lie
u\cap\lie k,V_{(\ka,K)}^{\mu+\rho_c}),S^n(\lie u\cap\lie
p)\otimes V_{(\lie m,M)}^{\lambda_1}).
$$
By Corollary 3.8 of \cite{KV}, as $M$-modules,
$$
H_j(\lie
u\cap\lie k,V_{(\ka,K)}^{\mu+\rho_c})\simeq H^{s-j}(\lie
u\cap\lie k,V_{(\ka,K)}^{\mu+\rho_c}\otimes\Lambda^{{\text{top}}}\lie
u\cap\lie k).
$$
Since the action of $\lie u\cap\lie k$ on $\Lambda^{{\text{top}}}\lie
u\cap\ka$ is
trivial we can write
$$
H_j(\lie
u\cap\lie k,V_{(\ka,K)}^{\mu+\rho_c})\simeq H^{s-j}(\lie
u\cap\lie k,V_{(\ka,K)}^{\mu+\rho_c})\otimes\Lambda^{{\text{top}}}\lie
u\cap\lie k,
$$
hence, using the fact that $V_{(\lie
m,M)}^{\lambda_1}\otimes(\Lambda^{{\text{top}}}\lie u\cap
\lie k)^*=V_{(\lie m,M)}^{\mu_\lambda+\rho_{\lie m}}$, we find that
$$
m_\mu=
\sum_{j=0}^s(-1)^{s-j}\sum_{n=0}^\infty\dim\hbox{\rm Hom}_M(H^{s-j}(\lie
u\cap\lie k,V_{(\ka,K)}^{\mu+\rho_c}),S^n(\lie u\cap\lie
p)\otimes V_{(\lie m,M)}^{\mu_\lambda+\rho_{\lie m}}).
$$

Set $W^1=\{w\in W_\tau\mid N(w)\subset\Delta(\lie u\cap\ka)\}$.
Applying Kostant's Theorem we find that
$$
m_\mu=\sum_{j=0}^s(-1)^{s-j}\sum_{n=0}^\infty\!\!\sum_{ w\in
W^1\atop\ell(w)=s-j}\!\!\!\!\!\dim\hbox{\rm Hom}_M(V_{(\lie
m,M)}^{w(\mu+\rho_c)-\rho_c+\rho_{\lie m}},S^n(\lie u\cap\lie
p)\otimes V_{(\lie m,M)}^{\mu_\lambda+\rho_{\lie m}}).
$$
By Lemma~2.7, if
$$
\dim\hbox{\rm Hom}_M(V_{(\lie
m,M)}^{w(\mu+\rho_c)-\rho_c+\rho_{\lie m}},S^n(\lie u\cap\lie
p)\otimes V_{(\lie m,M)}^{\mu_\lambda+\rho_{\lie m}})\ne0,
$$
then $w(\mu+\rho_c)-\rho_c$ is dominant for $\Delta_\tau^+$, hence
$w(\mu+\rho_c)$ is dominant and regular. Since $\mu+\rho_c$ is dominant and
regular, we find that $w=1$ and the above sum  reduces to (5).
\endemo

Theorem~2.8 says that to compute the $K$-spectrum of
$V^\lambda_{(\g,K)}$ it is enough to
compute the $M$-spectrum of $S(\lie u\cap\lie p)$. The point is the
fact that the
$M$-spectrum of $S(\lie u\cap\lie p)$ is somewhat easy to compute.
Indeed  suppose
that $\lambda\in P_{\text{reg}}(z)$ with $z\in \widetilde
Z^{\text{cmpt}}_\tau$ and
let $\tau'$ be as in the proof of Lemma~2.7. We classify $\tau'$
according to Table~I
below. If $2\tau$ is of type 1 then $\lie m=\lie k$ and $\lie u=\lie
u\cap\lie p$. This
means that $V^\lambda_{(\lie g,K)}$ is an holomorphic discrete series and the
$K$-structure of $S(\lie u)$ is very well known.

If $\tau'$ is of type 2, then $V^\lambda_{(\g,K)}$ is a Borel-de
Siebenthal discrete
series and the $M$-spectrum of $S(\lie u\cap\lie p)$ is given in
\cite{M-F} for the
classical cases and $F_4$.

We now discuss the case of $\tau'$ of type 3, that is
$\tau'=\o_i+\o_j$ with $i\ne j$ and
$i,j\in J$. Set
$$
\D_i=\{\alpha\in\Dp\mid (\o_i,\alpha)=1,\ (\o_j,\alpha)=0\},
$$
$$
\D_j=\{\alpha\in\Dp\mid (\o_j,\alpha)=1,\ (\o_i,\alpha)=0\}.
$$
and $\lie u_i$ (resp. $\lie u_j$) be the subspace of $\lie u\cap\lie
p$ such that
$\D(\lie u_i)=v_z^{-1}\Delta_i$ (resp. $\D(\lie u_j)=v_z^{-1}\D_j$).
Then, as a $M$-module,
$\lie u\cap\lie p=\lie u_i\oplus\lie u_j$, hence
$$
S^n(\lie u\cap\lie p)=\sum_{h+k=n}S^h(\lie u_i)\otimes S^k(\lie u_j).
$$

Set $\lie u_i^-=\sum_{\alpha\in v_z^{-1}\D_i}\g_{-\alpha}$, $\lie
u_j^-=\sum_{\alpha\in v_z^{-1}\D_j}\g_{-\alpha}$, and
$$
\g_i=\lie u^-_i\oplus\lie m\oplus\lie u_i\quad\g_j=\lie
u_j^-\oplus\lie m\oplus\lie u_j.
$$
We notice that $\lie u_i$ and $\lie u_j$ are abelian in $\g_i$ and
$\g_j$ respectively,
thus we can apply \cite{S1} to compute the $M$-spectrum of $S(\lie
u_i)$ and $S(\lie u_j)$.

\def\ep{\epsilon}
As an example we compute the $M$-spectrum of $S(\lie u\cap\lie p)$
for $G$ of type $E_6$
and $\tau'=\o_1+\o_6$. We use the notations of \cite{Bou}.
In this case
$$
\D_1=\{\half(\ep_8-\ep_7-\ep_6-\ep_5)+\half\sum_{i=1}^4(-1)^{\nu(i)}\ep_i\text{
with $\sum_{i=1}^4
\nu(i)$ odd}\},
$$
$$
\D_6=\{\pm\ep_i+\ep_5\mid 1\le i\le 4\}.
$$
If we set $\D_0=\{\pm\ep_i\pm\ep_j\mid 2\le i<j\le4\}$ and  $\lie m'=
\lie h\oplus\sum_{\alpha\in\D_{0}}\g_{\alpha}
$, then $\lie m=v_{z}^{-1}\lie m'$.

We set
$$
\align
\eta_{1}&=\half(\ep_8-\ep_7-\ep_6-\ep_5)\\
\eta_{2}&=\half(-\ep_{1}+\ep_{2}+\ep_{3}+\ep_{4})\\
\eta_{3}&=\half(\ep_{1}-\ep_{2}+\ep_{3}+\ep_{4})\\
\eta_{4}&=\half(\ep_{1}+\ep_{2}-\ep_{3}+\ep_{4})\\
\eta_{5}&=\half(\ep_{1}+\ep_{2}+\ep_{3}-\ep_{4}),
\endalign
$$
then we see that
$\{v_{z}^{-1}(\eta_{i}\pm\eta_{j})\mid 1\le i<j\le
5\}$ is the set of  positive
roots for $(\g_{1},\lie h)$ contained in $v_{z}^{-1}\Dp$ and that
$v_z^{-1}(\eta_{1}-\eta_{2})$ is the unique simple noncompact root. Applying
the results of \cite{S1}, we see that  the highest weights of
the $M$-types of $S(\lie u_{1})$ are given by
$\mu(h_{1},h_{2})=v_{z}^{-1}(h_{1}(\eta_{1}+\eta_{2})+h_{2}(\eta_{1}-\eta_{2}))$,
with $h_{1}\ge h_{2}$.

Similarly the positive system for $(\g_{6},\lie h)$ contained in
$v_{z}^{-1}\Dp$ is $\{v_{z}^{-1}(\pm\ep_{i}+\ep_{j})\mid 1\le i<j\le
5\}$ and $v_z^{-1}(\ep_{5}-\ep_{4})$ is the unique noncompact simple root. It
follows that the highest weights of
the $M$-types of $S(\lie u_{6})$ are given by
$\nu(k_{1},k_{2})=v_{z}^{-1}(k_{1}(\ep_{4}+\ep_{5})+k_{2}(\ep_{5}-\ep_{4}))$,
with $k_{1}\ge k_{2}$.
Putting all together we obtain that
$$
S^{n}(\lie u\cap\lie p)=\sum_{h_{1}+h_{2}+k_{1}+k_{2}=n\atop
h_{1}\ge h_{2},\ k_{1}\ge
k_{2}}F_{M}(\mu(h_{1},h_{2}))\otimes F_{M}(\nu(
k_{1},k_{2})),
$$
where $F_{M}(\mu)=V_{(\lie m, M)}^{\mu+\rho_{\lie m}}$ denotes the $M$-type of
highest weight $\mu$.
\medskip

\heading \S3 Nilradical and special abelian ideals\endheading
\medskip
If $\tau \in \h_{\Bbb R}$
    we set
$\frak{q}_\tau=\frak{m}_\tau\oplus\frak{n}_\tau$,
where $\frak{m}_\tau=\h\oplus\sum\limits_{(\alpha,\tau)=0}\g_\a$ is a
Levi subalgebra and
$\n_\t=\sum\limits_{\a\in\Dp:(\a,\t)>0}\g_\a$ is its nilradical.
\smallskip\noindent
\remark{Remark}
If $\i$ is abelian, then $\n_{\langle\i\rangle}$ is an $ad$-nilpotent
ideal of $\bb$; moreover
$\i\subset \n_{\langle\i\rangle}$. The first assertion follows from
the fact that $\langle\i\rangle$ is a dominant
weight (cf. \cite{Ko1, Prop. 6}). For the second claim, remark that,
since $\a+\b$ is not a root  for any $\a,\b\in\p_\i$,
we have $(\a,\b)\geq 0$; therefore, for $\a\in\p_\i$, we have
$(\a,\langle\i\rangle)=(\a,\a)+\sum\limits_{\b\in\p_\i\setminus\{\a\}}\,(\a,\b)>0$.
\endremark

In the following Lemma we recall two well known facts about roots.

\proclaim{Lemma 3.1} (1).
If $\a, \b\in \Dp$ and $\a<\b$, then there exists $\b_1,\dots,
\b_h\in \Dp$ such that
$\a+\sum\limits_{i=1}^j\b_i\in \Dp$ for $j=1, \dots , h$, and
$\a+\sum\limits_{i=1}^h\b_i=\b$.
\smallskip\noindent
(2).
If $\a+\b, \gg, \a+\b+\gg\in \D$, then at least on of  $\a+\gg$,
$\b+\gg$ is a root.
\endproclaim
\demo{Proof}
 To prove (1), fix $\beta\in \Delta^+$. We will show that for each $\alpha\in\Delta^+$ such that
$\alpha<\beta$, there exist $\beta_1, \dots, \beta_h\in \Delta^+$ such that 
$$
\beta=\alpha+\beta_1+\dots+\beta_h
$$
and $\alpha+\beta_1+\dots+\beta_i\in\Delta^+$ for all $i=1,\dots,h$.
We shall prove the statement by induction on $ht(\b-\a)$.
If $ht(\b-\a)=1$ the result is obvious.
Suppose that $ht(\b-\a)>1$.
If $(\a, \b)>0$, then $\b-\a$ is a root and we are done, so we
assume $(\a, \b) \leq 0$. By assumption, there exist
$\eta_1, \dots, \eta_k\in \Dp$ such that $\b=\a+\eta_1+\cdots+\eta_k$
and we obtain that, for at least one of the $\eta_i$,
$(\a, \eta_i)< 0$. We may assume $(\a, \eta_1)< 0$, so that
$\a+\eta_1\in \Dp$. If $\a+\eta_1=\b$ we are done, otherwise $\a+\eta_1<\b$ and  we
can apply the induction hypothesis.

 Assertion (2) is
an immediate consequence of the Jacobi identity.
\endemo
\smallskip
Recall that  $X=\overline C_2\cap P^\vee$. Let $\tau\in X$ and 
set $\i_{(\tau)}=\bigoplus\limits_{\a\in \D_\tau^2} \g_\a$. Then 
$\i_{(\tau)}$
is an abelian ideal of $\bb$ that should not be confused with
$\i_{\tau}$ in case $\tau\in X\cap\widetilde{Z}_{ab}$.

\proclaim{Lemma 3.2}
If $\tau\in X$ and  $\D_\tau^2\not=\emptyset$ then
\roster
\item
for all $\a\in \D_\tau^1$ there exists $\b\in\D_\tau^1$ such that
$\a+\b\in \D_\tau^2$.
\item
$
\n_{\langle\i_{(\tau)}\rangle}=\n_\tau=\sum\limits_{\a\in\D_\tau^1\cup\D_\tau^2}\g_\a.$
\endroster
\endproclaim

\demo{Proof}
\noindent

(1) Assume $\a\in \D^1_\tau$. Consider $\D(\a)=\{\b\in\D^2_\tau\mid \b>\a\}$.
Since
$(\tau,
\th)=2$ we have
$\a\not=\th$, hence $\D(\a)\ne\emptyset$. Pick $\b$ minimal in $\D(\a)$.

By Lemma 3.1 (1), we can find
$\b_1, \dots, \b_k\in \Dp$ such that
$\gamma_j=\a+\sum\limits_{i=1}^j\b_i$ is a root for
$j=1, \dots, k$ and $\gamma_k=\b$.   Choose among such expressions one such that
$k$ is minimal. By the choice of $\beta$, 
$(\b_k, \tau)=1$. We shall prove the result by  showing that $k=1$.
  If $k>1$, 
we have $\gamma_k=\gamma_{k-2}+\b_{k-1}+\b_k$, where we set $\gamma_0=\a$.
By Lemma 3.1 (2),either $\gamma_{k-2}+\b_k\in\Dp$ (hence $\gamma_{k-2}+\b_k\in\D^2_\t$), but this is
not possible by the minimality of $\beta$, or $\b_{k-1}+\b_k\in
\Dp$, and indeed $\b_{k-1}+\b_k\in \D_\tau^1$.
But in this latter case, setting $\b'_{k-1}=\b_{k-1}+\b_k$, we see that
$\b=\a+\b_1+\dots+\beta_{k-2}+\b'_{k-1}$ contradicting the minimality of $k$.
\par\noindent
(2) Set $\i=\i_{(\tau)}$. We first prove the inclusion 
$\D_\tau^0\subset \D_{\langle\i\rangle}^0$.
since $\D_\tau^0$ is standard parabolic subsystem of $\D$, it
suffices to prove that $\D_\tau^0\cap
\Pi\subset\D_{\langle\i\rangle}^0$.
If $\a\in \D_\tau^0\cap \Pi$ and $\b\in \D^2_\tau$, then clearly
$s_\a(\b)\in \D^2_\tau$.
This implies that
$s_\a\(\sum\limits_{\b\in\D^2_\tau}\b\)=\sum\limits_{\b\in\D^2_\tau}\b$,
hence
that $\(\a, \sum\limits_{\b\in\D^2_\tau}\b\)=0$, or $(\a,\langle\i\rangle)=0$.
It remains to prove the reverse inclusion:
we assume that $(\a, \tau)>0$ and prove that $(\a,\langle\i\rangle)>0$.
If $\a\in \D^2_\tau$, then by definition  $\a\in\p_\i$, hence, by the
above Remark, $(\a,\langle\i\rangle)>0$.
Thus it suffices to prove that if $\a\in \D^1_\tau$, then
$(\a,\langle\i\rangle)>0$.
If $\a\in \D^1_\tau$ and $\gamma\in\D^2_\tau$ then necessarily
$(\a,\gamma)\geq 0$: otherwise
$\a+\gamma$ is a root and $(\a+\gamma,\tau)=1+2=3$, a contradiction.
It remains to prove that for at least one $\gg\in \D^2_\tau$ we have
$(\a,\gamma)> 0$. By (1)
we can find $\b\in \D_\tau^1$ such that $\a+\b\in\D_\tau^2$; we shall
prove that then $(\a, \a+\b)>0$.
Assume by contradiction $(\a, \a+\b)\leq 0$. Then in particular
$(\a, \b)<0$, so that
$(s_\a(\b), \tau)=(\b-(\b, \a^\vee)\a, \tau)\geq 2$. It would follow that
$(\b-(\b, \a^\vee)\a, \tau)= 2$, hence that $(\b, \a^\vee)=-1$ and
therefore  that
$(\a^\vee, \a+\b)=1$: this is impossible since $(\a^\vee, \a+\b)$
differs from $(\a, \a+\b)$
by a positive factor.
\endemo

\proclaim{Theorem 3.3}\cite{Ko2, Theorem 4.4}
Given $\i\in\I_{ab}$, the following conditions are equivalent:
\roster
\item $[\n_{\langle\i\rangle},
\n_{\langle\i\rangle}]\subset\i=\cent \n_{\langle\i\rangle}$.
\item $\p_\i=\D^2_\tau$ for some $\tau \in X$.
\endroster
\endproclaim

\demo{Proof}
($1 \Rightarrow 2$) We already remarked, as a general fact, that
$\i\subset \n_{\langle\i\rangle}$.
We assume that condition (1) holds and for $\a\in \Dp$ we define
$$\tau(\a)=\cases 2 \quad \text{ if  }\quad \g_\a\subset\i \\
                     1 \quad \text{ if  }\quad  \g_\a\subset
\n_{\langle\i\rangle}\setminus \i \\
                     0 \quad \text{ if  }\quad
\g_\a\not\subset\n_{\langle\i\rangle}. \\
            \endcases$$
We shall prove that if $\a, \b, \a+\b\in \Dp$, then
$\tau(\a+\b)=\tau(\a)+\tau(\b)$.
We first verify that if $\g_\a\not\subset\n_{\langle\i\rangle}$,
$\g_\b\subset\n_{\langle\i\rangle}$,
and $\a+\b\in
\Dp$, then $\g_{\a+\b}\subset \n_{\langle\i\rangle}$, and
$\g_{\a+\b}\subset \i$ if and only if
$\g_\b\subset \i$.
The first condition is immediate since $\n_{\langle\i\rangle}$ is an
ideal of $\frak q_{\langle\i\rangle}$.
Since $\i=\cent\n_{\langle \i\rangle}$,  we have in particular that $\i$ is an ideal of
$\frak q_{\langle\i\rangle}$, too;  hence,  if
$\g_\b\subset\i$, $\g_{\a+\b}\subset\i$, too. Similarly, if
$\g_\b\subset
\n_{\langle\i\rangle}\setminus \i$, then $\g_{\a+\b}\subset
\n_{\langle\i\rangle}\setminus \i$, otherwise we should
obtain that $\g_\b=[\g_{-\a}, \g_{\a+\b}]\subset\i$.
The next case to consider is when $\g_\a,\g_\b$ are both contained in $\n_{\langle \i\rangle}\setminus \i$; if $\a+\b\in \D$,
the first equality in (1) implies $\g_{\a+\b}\subset\i$, as desired.
Finally, if $\g_\a\subset \i=\cent \n_{\langle\i\rangle}$ and
$\g_\b\subset \n_{\langle\i\rangle}$ we have that
$\a+\b\not\in \D$, so we can conclude that $\tau$ is additive and less or equal
than $2$ on $\Dp$.
This fact implies that $\tau$ can be extended to a linear functional
on $\h^*_\real$ which clearly corresponds to an element in $X$; if we
still denote by $\tau$ this
element, we obtain by definition that $\p_\i=\D_\tau^2$.
\par
($2 \Rightarrow 1$)
By Lemma 3.2, (2) we have
$\n_{\langle\i\rangle}=\bigoplus\limits_{\a\in\D_\tau^1\cup\D_\tau^2}\g_\a$.
Hence it
is
clear that $[\n_{\langle\i\rangle},
\n_{\langle\i\rangle}]\subset\i$; it is also clear
that
$\i \subset \cent \n_{\langle\i\rangle}$. Item (1) of Lemma 3.2 shows that no
element in
$\bigoplus\limits_{\a\in\D_\tau^1}\g_\a$ belongs to $\cent
\n_{\langle\i\rangle}$, so
indeed
$\i =\cent \n_{\langle\i\rangle}$.
\endemo

\proclaim{Definition} An abelian ideal of $\bb$ is said special if it satisfies
one of the conditions of Theorem 3.3.
An abelian ideal is  nilradical  if it is
the nilradical of a parabolic subalgebra of $\g$.
\endproclaim
\medskip Now we determine the structure of both nilradical and
special abelian ideals. In particular
we provide a proof of Theorem 4.9 of \cite{Ko2}.\par
Set  ${\Cal M}=\{\o_i\mid i\in J\}\cup\{0\}$. It is well known 
that ${\Cal M}$ is a set of
representatives for $P^\vee/Q^\vee$.
\medskip\proclaim{Proposition 3.4} An abelian ideal $\i$ of $\bb$ is
nilradical if and
only if there exists exists $\omega\in\Cal M$ such that
    $\i=\frak{n}_{\omega}$. In particular nilradical abelian
ideals are in bijection with $Cent(G)$.
\endproclaim
\demo{Proof}
The only nontrivial thing to prove is the fact that
$0\ne\frak{n}_\tau=\bigoplus\limits_{\a\in\Dp:(\a,\tau)>0}\g_\a$ is
abelian only if we can choose
$\tau=\omega^\vee_i$ for some
$i\in J$.\par
Let $\alpha\in\Dp$ be such that $(\alpha,\tau)>0$ and such that
$\alpha=\alpha'+\alpha''$ with $\alpha', \alpha''\in\Dp$. Since $\n_\t$
is abelian,
$(\alpha',\tau)$ and $(\alpha'',\tau)$ cannot be both nonzero. Suppose that
$(\alpha',\tau)=0$: then we have that $(\alpha'',\tau)>0$. Repeating this
argument we obtain that $\alpha=\alpha_i+\gamma$ where $\gamma$ is a sum of
positive roots
$\beta$ such that
$(\beta,\tau)=0$ and $\alpha_i$ is a simple root such that $(\alpha_i,\tau)>0$.
Apply now this observation to the highest root $\theta$:  then
$\theta=\alpha_i+\gamma$. This implies that $\tau=k\omega^\vee_i$
with $i\in J$,
hence
$\frak{n}_\tau=\frak{n}_{\omega^\vee_i}$.
\endemo

Recall that $\th=\sum_{i=1}^n m_i\a_i$.

\proclaim{Proposition 3.5} Suppose $\t\in X$.
Then we have the following five possibilities:
\par\noindent
(1) $\tau=2\o_i,\,i\in J$;
(2) $\tau=\o_i,\,m_i=2$;
(3) $\tau=\o_i+\o_j,\,i,j\in J$;
(4) $\tau=\o_i,\,i\in J$;
(5) $\tau=0.$
Then $\D_\tau^2\not=\emptyset$, whence $\i_{(\tau)}\not=0$,
exactly in cases (1), (2), (3).\par In case (1) the special abelian
ideal $\i_{(\tau)}$ is nilradical (and
viceversa);  in  cases (2) and (3) $\i_{(\tau)}$ is not  nilradical  and relations
$\i_{(\tau)}\subsetneq \n_\tau$, $\i_{(\tau)}=[\n_\tau,\n_\tau]$ hold.
\endproclaim

\demo{Proof}
Since $\t\in \overline C_2$ we have $\t=\sum\limits_{i=1}^n
\varepsilon_i o_i$ with
$\varepsilon_i\geq 0$  and $0<\sum\limits_{i=1}^ n \varepsilon_i\leq 2$
  where
$o_i=\o_i/m_i,\,i=1,\dots, n$.
Since, moreover, $\t\in P^\vee=\bigoplus\limits_{i=1}^n\Bbb Z \o_i$,
the only possibilities
for $\t$ are the five ones listed above. It is clear that
$\D^2_\tau=\emptyset$ if and only if
we are in cases (4) or (5).
Suppose  that  $\tau=2\o_i,\,i\in J$:  then $\D^1_\tau=\emptyset$ and
part (2) of  Lemma 3.2
implies that $\i_{(\tau)}=\frak n_{\langle\i_{(\tau)}\rangle}$, hence 
$\i_{(\tau)}$
is nilradical.
It is clear that in cases (2) and (3) $\D^1_\tau\not=\emptyset$,
whence, by Lemma 3.2, (2),
$\i_{(\tau)}\subsetneq \n_{\langle\i_{(\tau)}\rangle}$; moreover, by 
Theorem 3.3,
$\i_{(\tau)}\supset 
[\n_{\langle\i_{(\tau)}\rangle},\n_{\langle\i_{(\tau)}\rangle}]$.
So we are left  with proving the reverse inclusion. It suffices to
prove that any $\b\in\D^2_\tau$
splits as the sum of two roots in $\D^1_\tau$. We note that in cases
(2) and (3)
if $\b\in\D^2_\tau$ then $\b$ is not simple, so $\b=\b_1+\b_2$ with
$\b_1, \b_2\in\Dp,\,\b_1, \b_2<\b$. If
$\b_1$ and $\b_2$ belong to $\D^1_\tau$ we are done; otherwise one of
them, say for
example $\b_1$, belongs to $\D^2_\tau$. By induction on the height of
$\b$ we can assume
that there exist $\gg_1, \gg_2\in \D^1_\tau$ such that
$\b_1=\gg_1+\gg_2$. Then by Lemma 3.1,
(2) at least one among $\gg_1+\b_2$, $\gg_2+\b_2$ is a root which
clearly belongs to
$\D_\tau^1$ and we get our claim. From the previous Proposition it
follows that  $\i_{(\tau)}$ is not  nilradical.
\endemo

Putting together the results of Lemma 3.2 and Theorem 3.3 we obtain
the following result.

\proclaim{Corollary 3.6}
The special abelian ideals of $\bb$ are exactly the
$\i_{(\tau)}$, with $\tau\in X$.
For any $\tau\in X$  we have $[\n_\tau,
\n_\tau]\subset \i_{(\tau)}=\cent \n_\tau$.
\endproclaim

\remark{\bf Remark} Lemma~3.2 allows easily to classify and
enumerate the special abelian ideals by looking at the coefficients of the
highest root. Indeed, if $\i$ is a nonzero special abelian ideal and
$\i=\i_{(\tau)}=\i_{(\tau')}$, then part (2) of Lemma~3.2 says that $\lie 
n_\tau=\lie
n_{\tau'}$, so
$\Delta^0_\tau=\Delta^0_{\tau'}$. Since we are assuming that
$\D^2_\tau=\D^2_{\tau'}$, we can conclude that $\D^1_\tau=\D^1_{\tau'}$ also.
This implies that $\tau=\tau'$.
  With reference to the five possibilities for $\tau$ listed Proposition 3.5,
the list of the  $\tau\in X$ such that $\i_{(\tau)}\ne0$ is  given in 
Table I; we label
Dynkin diagrams as in \cite{Bou}.
\endremark

\par\newpage

\heading Table I\endheading
$$
\alignat4
&\text{$\D$}\qquad\qquad&&\text{$\t$ of type 1}\qquad\quad&&\text{$\t$ of type
2}\qquad\quad&&\text{$\t$ of type
3}\\  \\
&A_n&&2\o_1,\ldots,2\o_n&& &&\o_i+\o_j,\ 1\leq i<j\leq n\\
&B_n&&2\o_1 &&\o_2,\ldots,\o_n\\
&C_n&&2\o_n &&\o_1,\ldots,\o_{n-1}\\
&D_n&&2\o_1,2\o_{n-1},2\o_{n} &&\o_2,\ldots,\o_{n-2}&&
\o_1+\o_{n-1},\,\o_1+\o_{n},\,\o_{n-1}+\o_{n}\\
&E_6&&2\o_1,2\o_6&&\o_2,\o_3,\o_5&&\o_1+\o_6\\
&E_7&&2\o_7&&\o_1,\o_2,\o_6\\
&E_8&& &&\o_1,\o_8\\
&F_4 && &&\o_1,\o_4\\
&G_2&& &&\o_2\\
\endalignat
$$

For $\a\in\Dp$ we set $V_\a=\{\b\in\Dp\mid \b\geq \a\}$.
  Let $\i_{(\tau)},\ \t\in X$, be a special abelian ideal. Then
  if
$\t=2\o_i,$ $i\in J$ or
$\t=\o_i,\,m_i=2$ we have
$\n_{\tau}=\bigoplus\limits_{\b\in V_{\a_i}}\g_\b$, whereas if 
$\t=\o_i+\o_j,\,i,j\in J$ then
$\n_{\tau}=\bigoplus\limits_{\b\in V_{\a_i}\cup V_{\a_j}}\g_\b$.
\bigskip\bigskip
We want to restate the previous results in terms of the machinery 
introduced in section 1.
If $\tau\in X$, set $\frak k=\frak k_\tau$ and
$\frak{p}=\frak{p}_\tau$. Set also $\frak{b}_{\frak k}=\frak{b}\cap
\frak k$, $\frak{b}_{\frak{p}}=\frak{b}\cap\frak{p}$ and
$\i(\p)=\bigoplus\limits_{\a\in\p}\g_\a$, for $\p\subset \Dp$.
Recall that $\Phi^1_z=v_z^{-1}(\Delta^{1}_z)$,
$\Phi^2_z=-v_z^{-1}(\Delta^{-2}_z)$, so that $\i_z=\i(\Phi^1_z)\oplus 
\i(\Phi^2_z)$. Moreover
we have 
$\bb_\ka=\h\oplus\left(\bigoplus\limits_{\a\in\D_\tau^0\cup\D_\tau^2}\g_\a\right),\,
\bb_{\frak p}=\bigoplus\limits_{\a\in\D_\tau^1}\g_\a$ and
$
\Delta^2_\tau=-v_z^{-1}(\Delta^{-2}_z),\ 
\Delta_\tau^0=v_z^{-1}(\Delta^0_z),\ 
\Delta_\tau^1=v_z^{-1}(\Delta^1_z)\cup
-v_z^{-1}(\Delta_z^{-1}).$\par
 With this notation we have (see
\cite{Ko2, Prop.4.7, Theorem 4.9})

\proclaim{Proposition 3.7} Let $z\in\widetilde{Z}_{ab}$ and set
$\tau=dom(z)$. 
\roster
\item
The following relations are equivalent:\par\noindent
$i)\ \t=2\omega\text{ for some $\omega\in {\Cal M}$}$;\par\noindent
$ii)\ z=-2\omega\text{ for some $\omega\in {\Cal M}$}$;\par\noindent
$iii)\  n_\tau=1$;\par\noindent
$iv)\ \frak k_\tau=\g$.\par
\item $\i(\Phi^2_z)$ is a special abelian ideal.\hskip5cm\phantom{O}
\item $\i(\Phi^2_z)$ is  nilradical if and only if $\tau=2\omega$ for 
some $\omega
\in\Cal M$.
\item If $\i(\Phi^2_z)$ is not  nilradical  then
$\i(\Phi^2_z)=[\frak{b}_{\frak{p}},\frak{b}_{\frak{p}}]$.
\item $\i(\Phi^1_z)=\frak{b}_z\cap\frak{b}_{\frak{p}}$.
\item $\i(\Phi^1_z)=\i_z\cap \frak{p}$ and
$\i(\Phi^2_z)=\i_z\cap\frak k$.
\endroster\endproclaim
\demo{Proof}
$i)\Rightarrow ii)$: if $\alpha\in \Delta^+$ then
$\alpha(z)=v_z^{-1}(\alpha)(\tau)\in\{0,-2\}$, hence $z/2$ is minuscule.
$ii)\Rightarrow iii)$: $\Delta_\tau=\Delta$, hence $n_\tau=1$.
$iii)\Rightarrow iv)$:  $W_\tau=W$ implies that $\Delta=\Delta_\tau$ 
and in turn that  $\frak k_\tau=\g$.
$iv)\Rightarrow i)$:  if $\Delta_\tau=\Delta$ then, if $\alpha\in\Delta^+$,
$(\alpha,\tau)\in\{0,2\}$ hence $\tau/2$ is minuscule. At this point part (1)
is completely proved. The other parts follow immediately combining 
the relations
listed just before the statement of this Proposition with Lemma 3.2 
and  Proposition 3.4.
\endemo

\medskip
\proclaim{Proposition 3.8}
Let $\tau\in X$, $\ss$ be an $\h$-submodule of $\bb_\pp$,
$\D(\ss)=\{\a\mid\g_\a\subset \ss\}$, and set
$\i=\i(\D^2_\tau)\oplus \ss$. Then the following facts are equivalent.
\roster
\item $\ss$ is a $\bb_\ka$-submodule of $\bb_\pp$;
\item $\i$ is an ad-nilpotent ideal of $\bb$ included in $\n_\tau$;
\item $\i$ is an ad-nilpotent ideal of $\bb$ and $(\b, \tau)>0$ for 
all $\b$ in $\D(\ss)$.
\endroster
The same equivalences hold if we replace \lq\lq ad-nilpotent\rq\rq\ 
with \lq\lq abelian\rq\rq\
in (2) and (3) and consider abelian submodules in (1).
\endproclaim

\demo{Proof}
The equivalence of (2) and (3) is immediate from the definitions.
We prove the equivalence of (1) and (2).
By the definition of $\n_\tau$ and $\bb_\pp$ we have that in any case 
$\i\subset \n_\tau$;
By assumption $[\h, \ss]\subset \ss$; moreover, since $\ss\subset 
\i(\D_\tau^1)$,
we have $[\i(\D^2_\tau), \ss]=0$.
Moreover, $[\i(\D_\tau^1), \ss]\subset \i(\D^2_\tau)$, thus we obtain that
$\i$ is an (ad-nilpotent) ideal of $\bb$ if and only if 
$[\i(\D_0^\tau), \ss]\subset \ss$.
But by the definition of $\bb_\ka$, we also obtain that $\ss$ is a 
$\bb_\ka$-module
if and only $[\i(\D_0^\tau), \ss]\subset \ss$, so we get our claim.
It is clear that all the above arguments still work when restricting 
to abelian objects.
\endemo

\proclaim{Proposition 3.9}
Let $\i\in \I$, $\tau\in X$, and assume
$\i=\i(\D_\tau^2)\oplus \ss$, with $\ss \subset \i(\D_\tau^1)$.
Then there exists $z\in \widetilde Z$ such that $\i=\i_z$,  $\tau=dom(z)$,
and $\ss=\i(v_z^{-1}(\D_z^1))$.
\par\noindent
Conversely, if $\i=\i_z$ and $dom(z)\in X$, we have
$\i=\i\(\D_{dom(z)}^2\)\oplus \ss$, with $\ss \subset \i\(\D_{dom(z)}^1\)$.
\endproclaim

\demo{Proof}
Set $\Phi=\D(\ss)$, $\Phi^2=(\Phi+\Phi)\cap \Dp$, and
and $\Phi'=(\D_\tau^2\setminus \Phi^2)\cup (\D_\tau^1\setminus \Phi)$.
We first prove that there exists  $v\in W$ such that $\Phi'=N(v)$.
For this it suffices to prove that $\Phi'$ and its complement
in $\Dp$ are closed.
Since $\tau\in \overline C_2\cap P^\vee$, to prove  that $\Phi'$ is
closed it suffices to show that
if $\a,\b\in\D_\tau^1$ and $\a+\b\in\Phi^2$ then either $\a$ or $\b$
belongs to $\p$. Set $\gg=\a+\b$; since $\gg\in\Phi^2$,
there exist $\xi,\eta\in\Phi$ such that $\gg=\xi+\eta$. Now from
the expansion
$$0<(\gg,\gg)=(\a,\xi)+(\b,\xi)+(\a,\eta)+(\b,\eta)$$
we deduce that  one of
the summands in the right-hand side of the previous relation is positive.
Since the difference of two roots having positive scalar product is a
root, and since $\eta-\b=\a-\xi$ and $\eta-\a=\b-\xi$, we have that either
$\eta-\b\in\D$ or $\eta-\a\in\D$.
It suffices to consider the case
$\eta-\a\in\D$. Suppose $\a-\eta\in\Dp$ and remark that
$\a-\eta\in\D^0_\tau$;
since $\ss$ is a $\bb_{\frak k}$-module, we have
$\a=(\a-\eta)+\eta\in\p$ as desired. If instead
$\eta-\a\in\Dp$, we note that
$\eta-\a=\b-\xi\in\D$ and as above we deduce that $\b=(\b-\xi)+\xi\in\p$.\par
Next we consider
$\Dp\setminus \Phi'=\D_\tau^0\cup\Phi\cup \Phi^{2}$. Obviously,
$\D_\tau^0$ and $\Phi\cup \Phi^{2}$ are closed. If
$\a\in\D_\tau^0,\b\in\Phi$ and $\a+\b\in\D$
we have that $\a+\b\in\Phi$ since $\ss$ is a $\bb_{\frak k}$-module.
Finally, if $\a\in\D_\tau^0,\b\in\Phi^{2},
\b=\xi+\eta,\,\xi,\eta\in\Phi$, and
$\a+\b=\a+\xi+\eta\in\Dp$, we have by Lemma 3.1, (2), that either
$\a+\xi$ or  $\a+\eta$ is a root and hence belongs to
$\Phi$, so that $\a+\b\in \Phi^2$.
\par
So there exists  $v\in W$ such that  $\Phi'=N(v)$.
Set $z=v^{-1}(\tau)$. Then $N(v_z)=\{\a\in\Dp\mid
(\a,v^{-1}(\tau))<0\}=\{\a\in\Dp\mid (v(\a),\tau)<0\}=
\{\a\in\Dp\mid v(\a)<0\}=N(v^{-1})$, hence $v_z=v^{-1}$ and therefore $dom(z)=\tau$.
\par
In order to conclude we have to prove that $\i=\i_z$.
We have $\widetilde F (z)=t_\tau v$, hence, by Lemma C,
 $\i_z=
\i(\D_\tau^2)\oplus \i(\D_\tau^1\setminus N(v))$.
Since $\D(\ss)\subset \D_\tau^1$, and
$N(v)\cap \D_\tau^1=\D_\tau^1\setminus \D(\ss)$, we have
$\D_\tau^1\setminus N(v)=\D(\ss)$, whence $\i=\i_z$.
Finally by a direct check we obtain that $\D_\tau^1\setminus 
N(v)=v(\D_z^1)$,
so $\ss=\i(v_z^{-1}(\D_z^1))$.
\par
The converse statement follows directly from Lemma C.
\endemo

\bigskip Recall that  $dom:\h_\Bbb R\to C_\infty$ is defined by 
$z\mapsto v^{-1}_z(z)$.
As already observed $dom(\widetilde{Z}_{ab})=X=
\overline{C}_2\cap P^\vee$, whereas,
by Proposition 1.3, we have $dom(\widetilde{Z})\subset\overline
C_{h}\cap P^\vee$.  We recall that
for $\tau\in X$ we set $\widetilde{Z}_\tau=dom^{-1}(\tau)\cap 
\widetilde{Z}_{ab}$;
we also set $\widehat{Z}_\tau=dom^{-1}(\tau)\cap \widetilde{Z}$.
Fix $\tau\in X$ and let $\Cal{S}_\tau$ denote the set of
$\frak{b}_{\frak k}$-submodules of $\frak{b}_{\frak{p}}$.

\medskip
\proclaim{Proposition 3.10}
The map $z\mapsto \i(v_z^{-1}(\D_z^1))$ establishes
a bijection between $\widehat{Z}_\tau$ and $\Cal{S}_\tau$. Moreover, this
map restricts to a bijection between $\widetilde{Z}_\tau$ and the 
abelian subspaces of
$\Cal{S}_\tau$. In particular, the number of
abelian  $\bb_\ka$-submodules of $\bb_\pp$ is $n_\tau$.
\endproclaim

\demo{Proof}
Since $\tau+Q^\vee=z+\Q$ for any $z\in \widehat Z_\tau$, the map
$z \mapsto \i_z$  is  injective on $\widehat{Z}_\tau$ (see 	Prop. 
1.4, (c)).
For any $z\in \widehat Z_\tau$ we have  $\widetilde F(z)=t_\tau v_z^{-1}$, 
hence, by Lemma C, $\Phi_{\i_z}=\D_\tau^2 \cup (\D_\tau^1\setminus
N(v_z^{-1}))$.
As in the proof of 3.9, we have $\D_\tau^1\setminus N(v_z^{-1})= 
v_z^{-1}(\D_z^1)$,
so that  $\i_z= \i(\D_\tau^2 )\oplus \i(v_z^{-1}(\D_z^1))$. In particular
$v_z^{-1}(\D_z^1)$ determines $\i_z$, and therefore
$z\mapsto \i(v_z^{-1}(\D_z^1))$ is one to one.
Moreover, by Proposition 3.8, we obtain that $\i(v_z^{-1}(\D_z^1))$ is a
$\bb_\ka$-submodule of $\bb_\pp$, so $z\mapsto \i(v_z^{-1}(\D_z^1))$ 
is a one to one map
from $\widehat{Z}_\tau$ to $\Cal{S}_\tau$. It remains to prove that this map
is onto. Take any $\ss\in \Cal S_\tau$; then, by Proposition 3.8,
$\D_\tau^2\oplus \ss$ is an ad-nilpotent ideal of $\bb$ included in $\n_\tau$.
Thus by Lemma 3.9 there exists $z\in \widetilde Z$ such that $dom(z)=\tau$
and $\ss=\i(v_z^{-1}(\D_z^1))$. So we have proved that $z\mapsto 
\i(v_z^{-1}(\D_z^1))$
is a bijection between $\widehat{Z}_\tau$ and $\Cal{S}_\tau$.
The final assertions are clear.
\endemo

\bigskip
Recall now the $Cent(G)$ action on $\widetilde Z$ introduced in 
Proposition 1.1,
explicitly given by the action of the subgroup $\Sigma$ of $\widetilde W$.
Let $\Sigma \cdot z$ denote the orbit of $z\in \widetilde Z$ under this action.
As before we denote by $\Cal M$ the set $\{\o_i\mid i\in J\}\cup\{0\}$ of
representatives for $P^\vee/Q^\vee$ and
for any $\omega\in \Cal M$ we set $V_\omega=\{\a\in \Dp \mid (\a, 
\omega)>0\}$. So
$V_{\o_i}=V_{\alpha_i}$ and $V_0=\emptyset$.
For $\i\in \I_{ab}$  we set
$$\align
&C_\i=\{\omega\in \Cal M\mid \p_\i\subset V_\omega\},\\
&C'_\i=dom(\Sigma\cdot z)\cap \Cal M,\\
\endalign
$$
where $z$ is any element of $\widetilde{Z}_{ab}$ such that $\i_{z}=\i$.
\medskip
\proclaim{Corollary 3.11}
\cite{Ko2, Theorem 5.1}
Let $\i\in \I_{ab}$ and let $z\in \Z_{ab}$ be such that $\i=\i_z$.
Then there are exactly $|Cent(G)|-|C_\i|$ decompositions
$$\i=\i(\D_\tau^2)+\ss,$$
with $\i(\D_\tau^2)$ nonzero special abelian ideal and $\ss\subset \i(\D_\tau^1)$.
\endproclaim

\demo{Proof}
By Proposition 1.4, for any $z'\in \Z$ we have $\i_{z'}=\i_z$ if and only if
$z'\in \Sigma\cdot z$.
Hence, by Proposition 3.9, we have $\i=\i(\D_\tau^2)+\ss,$ with $\tau\in X$ and
$\ss\subset \D_\tau^1$, if and only if $\tau=dom(z')$ for some $z'\in 
\Sigma\cdot z$. By the Remark following 3.6, the map $\t\to\D^2_\t$ is injective.
Therefore we have only to prove that $\D_{dom(z')}^2\not=\emptyset$ for
exactly $|Cent(G)|-|C_\i|$ elements $z'\in \Sigma\cdot z$.
We first notice that the restriction of $dom$ to a
$\Sigma$-orbit is injective. This follows from the fact that
distinct elements in $\Sigma\cdot z$ are distinct mod-$Q^\vee$,
as clear from the definition of $\Sigma$, while
$z'+ Q^\vee=dom(z')+Q^\vee$ for each $z'\in \Sigma\cdot z$.
In particular  we have $|dom(\Sigma\cdot z)|=|Cent(G)|$.
Now it suffices to prove that $\D_{dom(z')}^2=0$ for exactly $|C_\i|$ elements
$z'\in \Sigma\cdot z$.
It is clear that if $\tau\in X$, we have $\D_\tau^2=\emptyset$ if and only if
$\tau\in \Cal M$, therefore if we prove that $C'_\i=C_\i$ we are done.
Assume first $\omega\in C'_\i$. Then $\omega=dom(z')$ for some $z'$ such that
$\i=\i_{z'}$. Since $\widetilde F(z')=t_\omega v_{z'}^{-1}$ using  
Lemma C we obtain $\Phi_\i\subset V_\omega$.
 Conversely, assume $\omega\in C_\i$: then, clearly, $\D_\omega^2=\emptyset$ and 
  $\i\subset \i(\D_\omega^1)=\i(V_{\omega})$. Therefore, by Proposition 3.9,
there exists $z'\in \Sigma\cdot z$ such that $dom(z')=\omega$: this 
concludes the proof.
\endemo

 \enddocument

\bye